\documentclass[titlepage,final,11pt]{article}
\usepackage{theorem}
\usepackage{enumerate}
\usepackage{array}
\usepackage[reqno]{amsmath}
\usepackage{amssymb}
\usepackage{mathtools}
\usepackage{latexsym}
\usepackage{makeidx}
\usepackage{fancybox}
\input epsf.sty
\usepackage[dvips]{graphicx,color}
\usepackage{showlabels}
\usepackage{subcaption}

\theoremstyle{change}
\newtheorem{proclaim}{PROCLAIM}[section]
\newtheorem{theorem}[proclaim]{Theorem}
\newtheorem{definition}[proclaim]{Definition}

\newtheorem{proposition}[proclaim]{Proposition}
\newtheorem{corollary}[proclaim]{Corollary}
\newtheorem{example}[proclaim]{Example}

\numberwithin{equation}{section}

\outer\def\proclaim #1. #2\par{\medbreak \noindent{\bf#1.\enspace}{\sl#2}\par
  \ifdim\lastskip<\medskipamount
  \removelastskip\penalty55\medskip\fi}
\def\state #1. { \noindent{\bf#1.\enspace}}
\def\algo #1. { \noindent{\bf#1.\enspace}}

\DeclareMathOperator{\bdry}{bdry}
\DeclareMathOperator{\con}{con}

\DeclareMathOperator{\dist}{dist}

\DeclareMathOperator{\exs}{exs}
\DeclareMathOperator{\dom}{dom}
\DeclareMathOperator{\gph}{gph}
\DeclareMathOperator{\epi}{epi}

\DeclareMathOperator{\nt}{int}

\DeclareMathOperator{\prj}{prj}

\newcommand{\comp}{\,{\raise 1pt \hbox{$\scriptstyle\circ$}}\,}

\newcommand{\reals}{\mathbb{R}}
\newcommand{\Reals}{\overline{\mathbb{R}}}
\newcommand{\natnums}{{{\rm l} \kern -.13em {\rm N} }}
\newcommand{\nats}{\mathbb{N}}
\newcommand{\snats}{{I\kern -.29em N}}
\newcommand{\rats}{{Q\kern -.64em \raise 1pt \hbox{$\scriptstyle |$}\;\,}}
\newcommand{\srats}
	{{Q\kern -.56em \raise 1.2pt \hbox{$\scriptscriptstyle /$}\,}}
\newcommand{\ints}{Z\kern -.46em Z}
\newcommand{\ball}{\mathbb{B}}
\newcommand{\pluss}{\hskip1pt \raise1pt\vbox{\hrule width6pt \vskip1pt \hrule
                    width6pt} \kern-4pt{\lower1pt\hbox{\vrule height6pt
		    \kern1pt\vrule height6pt}}\hskip5pt}
\newcommand{\eop}
	{\hfill{$\vcenter{\hrule height1pt \hbox{\vrule width1pt height5pt
   	 \kern5pt \vrule width1pt} \hrule height1pt}$} \medskip}
\newcommand{\half}
	{{\raisebox{1pt}{$\frac{1}{2}$}}}
\newcommand{\third}
	{{\raisebox{1pt}{$\frac{1}{3}$}}}

\newcommand{\twothird}
	{{\raisebox{1pt}{$\frac{2}{3}$}}}

\newcommand{\setd}{{ d \kern -.15em l}}
\newcommand{\hatsetd}{ d \hat{\kern -.15em l }}

\renewcommand{\epsilon}{\varepsilon}
\renewcommand{\phi}{\varphi}

\newcommand{\lset}{\big\lbrace}

\newcommand{\rset}{\big\rbrace}

\newcommand{\tto}{\;{\lower 1pt \hbox{$\rightarrow$}}\kern -12pt
           \hbox{\raise 2.5pt \hbox{$\rightarrow$}}\;}
\newcommand{\overto}[1]{\,{\raise 0pt\hbox{$\rightarrow$}}\kern -9pt
     \hbox{\lower 3pt \hbox{$\scriptscriptstyle#1$}}\hskip6pt}
\newcommand{\underto}[1]{\,{\lower 1pt\hbox{$\rightarrow$}}\kern -9pt
     \hbox{\raise 4pt \hbox{$\,\scriptscriptstyle#1$}}\hskip7pt}
\newcommand{\bigoverto}[1]{{\raise 0pt\hbox{$\,\longrightarrow$}}\kern -16pt
     \hbox{\lower 3pt \hbox{$\scriptscriptstyle#1$}}\hskip4pt}
\newcommand{\bigunderto}[1]{\,{\lower 1pt\hbox{$\longrightarrow$}}\kern -16pt
     \hbox{\raise 4pt \hbox{$\,\scriptscriptstyle#1$}}\hskip6pt}
\newcommand{\bigbigto}[2]{\,{\raise 0pt\hbox{$\,\longrightarrow$}}\kern -16pt
     \hbox{\lower 3pt \hbox{$\scriptscriptstyle#2$}}\kern -10pt
     \hbox{\raise 4pt \hbox{$\,\scriptscriptstyle#1$}}\hskip7pt}
\newcommand{\downto}{{\raise 1pt \hbox{$\scriptscriptstyle \,\searrow\,$}}}
\newcommand{\upto}{{\raise 1pt \hbox{$\scriptscriptstyle \,\nearrow\,$}}}

\newcommand{\notimply}
	{\quad\hbox{$\Longrightarrow \kern -14pt {/}$}\hskip6pt\quad}

\newcommand{\lto}{\,{\lower 1pt\hbox{$\rightarrow$}}\kern -10pt
     \hbox{\raise 4pt \hbox{$\, \scriptstyle l$}}\hskip7pt}
\newcommand{\eto}{\,{\lower 1pt\hbox{$\rightarrow$}}\kern -10pt
     \hbox{\raise 4pt \hbox{$\, \scriptstyle e$}}\hskip7pt}
\newcommand{\hto}{\,{\lower 1pt\hbox{$\rightarrow$}}\kern -10pt
     \hbox{\raise 4pt \hbox{$\, \scriptstyle h$}}\hskip7pt}
\newcommand{\pto}{\,{\lower 1pt\hbox{$\rightarrow$}}\kern -10pt
     \hbox{\raise 4.5pt \hbox{$\, \scriptstyle p$}}\hskip7pt}
\newcommand{\cto}{\,{\lower 1pt\hbox{$\rightarrow$}}\kern -10pt
     \hbox{\raise 4pt \hbox{$\, \scriptstyle c$}}\hskip7pt}
\newcommand{\gto}{\,{\lower 1pt\hbox{$\rightarrow$}}\kern -10pt
     \hbox{\raise 4.5pt \hbox{$\, \scriptstyle g$}}\hskip7pt}
\newcommand{\sto}{\,{\lower 1pt\hbox{$\rightarrow$}}\kern -10pt
     \hbox{\raise 4pt \hbox{$\, \scriptstyle s$}}\hskip7pt}
\newcommand{\awto}{\,{\lower 1pt\hbox{$\rightarrow$}}\kern -15pt
     \hbox{\raise 4pt \hbox{$\, \scriptstyle aw$}}\hskip7pt}
\def\Nto{\,{\raise 1pt\hbox{$\rightarrow$}}\kern -13pt
     \hbox{\lower 3pt \hbox{$\, \scriptstyle N$}}\hskip7pt}
\def\Cto{\,{\raise 1pt\hbox{$\rightarrow$}}\kern -14pt
     \hbox{\lower 3pt \hbox{$\, \scriptstyle C$}}\hskip7pt}
\def\fto{\,{\raise 1pt\hbox{$\rightarrow$}}\kern -14pt
     \hbox{\lower 3pt \hbox{$\, \scriptstyle f$}}\hskip7pt}

\newcommand{\low}[1]{{\lower1pt \hbox{$\scriptstyle #1$}}}
\newcommand{\loww}[1]{{\lower2pt \hbox{$\scriptstyle #1$}}}
\newcommand{\high}[1]{{\raise1pt \hbox{$\scriptstyle #1$}}}

\newcommand{\cB}{{\cal B}}

\newcommand{\cN}{{\cal N}}

\newcommand{\nsum}{\mathop{\sum}\nolimits}

\newcommand{\nInnLim}{\mathop{\rm LimInn}\nolimits}
\newcommand{\nOutLim}{\mathop{\rm Lim\hspace{-0.01cm}Out}\nolimits}

\newcommand{\nliminf}{\mathop{\rm liminf}\nolimits}

\newcommand{\nlimsup}{\mathop{\rm limsup}\nolimits}
\newcommand{\ninf}{\mathop{\rm inf}\nolimits}
\newcommand{\nsup}{\mathop{\rm sup}\nolimits}
\newcommand{\nmax}{\mathop{\rm max}\nolimits}

\newcommand{\nnmin}{\mathop{\rm minimize}}

\newcommand{\nargmax}{\mathop{\rm argmax}\nolimits}
\newcommand{\nargmin}{\mathop{\rm argmin}\nolimits}

\newcommand{\bfxi}{\mbox{\boldmath $\xi$}}

\newcommand{\lwdy}[2]{\mathrel{\mathop
        {\raisebox{0.1ex}{\null$#1$}}{\hbox{\kern -1.0em
	{\raisebox{-0.8ex}{$\scriptstyle{\;\to #2}$}}}}}}
\newcommand{\lwwdy}[2]{\mathrel{\mathop
        {\raisebox{0.2ex}{\null$#1$}}{\hbox{\kern -1.0em
	{\raisebox{-1.1ex}{$\scriptstyle{\;\to #2}$}}}}}}
\newcommand{\slwdy}[2]{\scriptsize{{\mathrel{\mathop
        {\raisebox{0.1ex}{\null$#1$}}{\hbox{\kern -1.0em
	{\raisebox{-0.8ex}{$\scriptstyle{\;\to #2}$}}}}}}}}
\newcommand{\slwwdy}[2]{\scriptsize{{\mathrel{\mathop
        {\raisebox{0.2ex}{\null$#1$}}{\hbox{\kern -1.0em
	{\raisebox{-1.1ex}{$\scriptstyle{\;\to #2}$}}}}}}}}

\definecolor{lightgray}{gray}{0.75}
\definecolor{myred}{rgb}{0.55,0,0}
\definecolor{myblue}{rgb}{0,0,0.5} 
\definecolor{mygreen}{rgb}{0,0.5,0} 
\definecolor{purple}{rgb}{0.5,0,0.5} 
\definecolor{turq}{rgb}{0,0.805,0.816} 
\definecolor{maroon}{rgb}{0.51,0,0}
\definecolor{MAROON}{rgb}{0.51,0,0}
\definecolor{redor}{rgb}{0.78,0.078,0.078}
\definecolor{dgreen}{rgb}{0,0.3,0}

\newcommand{\Ex}{\mathbb{E}}

\newcommand{\grill}{{\scriptscriptstyle\#}}
\newcommand{\bcdot}{\,{\raise .2ex \hbox{$\centerdot$}}\,}

\newcommand{\bbA}{\mathbb{A}}

\newcommand{\bbP}{\mathbb{P}}

\begin{document}


\begin{center}
\begin{large}
{\bf Consistent Approximations in Composite Optimization}
\smallskip
\end{large}
\vglue 0.7truecm
\begin{tabular}{c}
  \begin{large} {\sl Johannes O. Royset } \end{large}\\
  Operations Research Department\\
  Naval Postgraduate School\\
  joroyset@nps.edu
\end{tabular}

\vskip 0.2truecm

\end{center}

\vskip 0.3truecm

\noindent {\bf Abstract}. \quad Approximations of optimization problems arise in computational procedures and sensitivity analysis. The resulting effect on solutions can be significant, with even small approximations of components of a problem translating into large errors in the solutions. We specify conditions under which approximations are well behaved in the sense of minimizers, stationary points, and level-sets and this leads to a framework of consistent approximations. The framework is developed for a broad class of composite problems, which are neither convex nor smooth. We demonstrate the framework using examples from stochastic optimization, neural-network based machine learning, distributionally robust optimization, penalty and augmented Lagrangian methods, interior-point methods, homotopy methods, smoothing methods, extended nonlinear programming, difference-of-convex programming, and multi-objective optimization. An enhanced proximal method illustrates the algorithmic possibilities. A quantitative analysis supplements the development by furnishing rates of convergence.\\


\halign{&\vtop{\parindent=0pt
   \hangindent2.5em\strut#\strut}\cr
                         \cr

{\bf Date}:\quad \ \today \cr}

\baselineskip=15pt

\section{Introduction}

A fundamental approach to optimization is to replace an actual problem by an approximating one, which is then solved using an existing algorithm. In sensitivity analysis, the actual problem of interest is also replaced by approximating ones for the purpose of identifying the effect of perturbations. These situations raise the questions: If some portion of an optimization problem is changed and the resulting problem is solved, would the obtained solution be a reasonable approximation of a solution of the actual problem? How large would the solution error be relative to the magnitude of the initial change? These questions are further complicated by the need for considering both optimal and stationary points in the nonconvex setting. Thus, we aspire to construct approximations that have minimizers  {\em and}\, stationary points near the corresponding points for the actual problem. We may even seek approximating objective functions and feasible sets that are near those of the actual problem in a broad sense.

In this paper, we provide a comprehensive framework for constructing and analyzing approximations of optimization problems with a composite structure. We provide sufficient conditions for {\em consistent approximations}, which guarantee that the approximating problems eventually become accurate relative to the actual problem in the sense of minimizers, stationary points, and level-sets. Examples from stochastic optimization, neural-network based machine learning, distributionally robust optimization, penalty and augmented Lagrangian methods, interior-point methods, homotopy methods, smoothing methods, extended nonlinear programming, difference-of-convex programming, and multi-objective optimization demonstrate the framework. The algorithmic possibilities are illustrated by a  proximal composite method for solving an array of nonconvex composite problems via a sequence of convex problems. A quantitative analysis supplements the development by furnishing rates of convergence.

We consider the broad class of composite optimization problems in the form
\begin{equation}\label{eqn:actualpr}
  \nnmin_{x\in \reals^n} ~\phi(x) = \iota_X(x) + h\big(F(x)\big),
\end{equation}
where the {\em objective function} $\phi:\reals^n\to \Reals = \reals \cup \{-\infty,\infty\}$ is defined in terms of a nonempty closed set $X\subset \reals^n$, often representing ``basic'' constraints such as bounds on the variables, and a locally Lipschitz continuous mapping $F:\reals^n\to \reals^m$ that models $m$ quantities of interest with $f_1, \dots, f_m$ being the {\em component functions}, i.e., $F(x) = (f_1(x), \dots, f_m(x))$. For any set $C$, $\iota_C(x) = 0$ if $x\in C$ and $\iota_C(x) = \infty$ otherwise. The quantities of interest are combined by a convex function $h:\reals^m \to \Reals$, which might assign an infinite penalty to certain values of $F(x)$. For example, $h = \iota_{\{0\}^m}$ sets $h(F(x)) = \infty$ whenever $F(x)$ deviates from the zero vector and thus encodes the equality constraint $F(x) = 0$. The actual problem \eqref{eqn:actualpr} captures classical nonlinear programming as well as many other problem formulations. In general, $\phi$ is nonconvex. The convexity assumption imposed on $h$ is not as restrictive as it might appear and helps us to leverage convexity properties when present.

The study of \eqref{eqn:actualpr} and similar problems goes at least back to \cite{Powell.78a,Powell.78b}. Recent efforts \cite{LewisWright.16,BurkeHoheisel.17,DrusvyatskiyLewis.18,DuchiRuan.18,CuiPangSen.18,DavisDrusvyatskiy.19,
DrusvyatskiyPaquette.19,BurkeHoheiselNguyen.21} suggest that the problem class might be among the most important ones in the nonconvex setting. Its structural properties are theoretically and computationally attractive, especially to address large-scale problems in machine learning and elsewhere.
Nevertheless, the actual problem may not be solvable directly and often needs approximations. The mapping $F$ may be approximated by a smooth or even an affine mapping. It might be specified by integrals, suprema, or other expressions that necessitate imprecise, numerical evaluation.  The function $h$ could be nonsmooth and extended real-valued, and then approximated by real-valued functions, piecewise affine or smooth, acting as penalties. The set $X$ might be approximated by a polyhedron or other simplifications. The approximations can be computationally motivated or introduced as part of sensitivity, stability, or error analysis.

Our approach is traced back to 1902 when P. Painlev\'{e} defined set-convergence. It is now well-known that set-convergence of epigraphs (i.e., epi-convergence) of approximating functions and set-convergence of graphs (i.e., graphical convergence) of approximating set-valued mappings furnish the desired convergence guarantees about the solutions of the corresponding minimization problems and generalized equations; see \cite[Chap. 4, 5, and 7]{VaAn}. If the set-valued mappings represent optimality conditions for the approximating problems, then graphical convergence translates into convergence of stationary points. We define approximating problems and their optimality conditions as being {\em consistent approximations} of the actual problem \eqref{eqn:actualpr}, paired with an optimality condition, when they exhibit such epi-convergence and graphical convergence.

In the convex setting, epi-convergence of approximating functions $\phi^\nu$ to the actual objective function $\phi$ suffices for the subgradient mappings $\partial \phi^\nu$ to graphically  converge to $\partial \phi$ by Attouch's theorem; see for example \cite[Thm. 12.35]{VaAn}. Thus, epi-convergence of the approximating functions ensures that the approximations are consistent. In the nonconvex case, one explicitly needs to impose conditions related to the optimality conditions to guarantee convergence of stationary points and this is recognized by the added requirement about graphical convergence of set-valued mappings in the definition of consistent approximations. If the approximating functions $\phi^\nu$ are smooth, their epi-convergence to $\phi$ does provide some relation between $\nabla \phi^\nu$ and $\partial \phi$ but in a less-than-ideal ``outer'' sense; cf. \cite[Lem. 3.4]{BurkeHoheisel.13} and \cite[Cor. 8.47]{VaAn}. A far-reaching extension of Attouch's theorem for nonconvex functions is provided by \cite{Poliquin.92}, which implies in the context of \eqref{eqn:actualpr} that epi-convergence of certain approximating functions $h^\nu$ to $h$ combined with a smooth $F$ and a constraint qualification suffice for graphical convergence in a local sense of the resulting subgradient mappings \cite{BurkeHoheisel.17}.

Approximations of \eqref{eqn:actualpr} may stem from the smoothing of $h$ and/or $F$. The specific function $\gamma\mapsto \max\{0,\gamma\}$ can be approximated by smooth functions via convolution  \cite{ChenMangasarian.95}. The approximations exhibit both epi-convergence as well as certain convergence of gradients to the subgradients of the max-function \cite{Chen.12}; the reference includes many examples of such smoothing. This approach is closely related to mollifiers, which may even approximate discontinuous functions \cite{ErmolievNorkinWets.95}. Further developments in these directions are furnished by \cite{BurkeHoheisel.13,BurkeHoheisel.17}, where in the context of \eqref{eqn:actualpr}, $h$ is approximated by a smooth function constructed using inf-convolution; see also \cite{BurkeHoheiselKanzow.13} for details about smoothing of finite max-functions. For approximations of a weakly convex expectation function caused by sample averages, \cite{DavisDrusvyatskiy.22} gives rates by which subdifferentials of the approximating functions graphically converge to the subdifferential of the expectation function.

In contrast to these efforts, we consider broad sets of approximations; $X$, $h$, and $F$ in \eqref{eqn:actualpr} may all be approximated and not necessarily smoothly. In sensitivity and error analysis, nonsmooth approximations are especially relevant as computational concerns tend to be secondary. We also deviate from the focus on the subgradient mapping of $\phi$ and that of its approximations, and instead express optimality conditions using (generalized) {\em multipliers}, which provide additional flexibility for absorbing inaccuracies. Thus, consistent approximations emerge as widely available under mild assumptions. As compared to \cite{DavisDrusvyatskiy.22} specifically, we consider a broader range of approximations beyond sample averages and also do not rely on weak convexity. While $\phi$ is weakly convex when $X=\reals^n$, $h$ is real-valued, convex, and Lipschitz continuous, and $F$ is continuously differentiable with Lipschitz continuous Jacobian, we emphasize situations with a nontrivial $X$, even nonconvex, an extended real-valued $h$, and nonsmooth $F$.

There is an extensive literature on general approximations and related stability analysis; see the monographs \cite{Polak.97,VaAn,BonnansShapiro.00,Mordukhovich.13b,CuiPang.21} as well as efforts based on metric regularity and calmness \cite{IoffeOutrata.08,Penot.10}, tilt-stability \cite{EberhardWenczel.12,LewisZhang.13,DrusvyatskiyLewis.13}, full-stability \cite{MordukhovichRockafellarSarabi.13}, convergence of abstract iteration schemes \cite{KlatteKrugerKummer.12}, and the truncated Hausdorff distance \cite{AttouchWets.91,AttouchWets.93a,AttouchWets.93b,Royset.18,Royset.20b}. The latter approach is closely related to the present paper, with direct relevance to our rate of convergence analysis.

Our choice of the term ``consistent approximations''  is motivated by E. Polak's concept in \cite[Chap. 3-4]{Polak.97}. There epi-convergence of approximating objective functions combined with an outer epi-limit of approximating optimality functions (after a sign change) are defined as consistent; see
\cite{Royset.12,RoysetPee.12,ForakerRoysetKaminer.16,PhelpsRoysetGong.16,RoysetWets.16c} for applications in stochastic optimization, semi-infinite programming, nonsmooth optimization, and optimal control.
Since optimality functions are essentially gap-functions of generalized equations representing optimality conditions, Polak's concept is, roughly, equivalent to weak consistency as we define it below. A related concept based on epigraphical nesting of directional derivatives is defined in \cite{HigleSen.95}. We deviate from these earlier developments by viewing optimality conditions as generalized equations defined by set-valued mappings and thus bypass the need for defining optimality functions. This brings into play well-developed calculus rules and computational procedures for generalized equations. The flexibility of this approach is demonstrated by a set of sufficient conditions for consistency (Theorem \ref{tConstentApproxComp}). An introductory treatment of consistent approximations under the assumption of a smooth mapping $F$ is, in parallel, provided by the textbook \cite{primer}.
While we limit the discussion to finite-dimensional problems, many of the concepts can be extended; see \cite{PhelpsRoysetGong.16, KouriSurowiec.20} for such efforts.

We continue in Section 2 by stating the approximating problems, defining consistency, and specifying sufficient conditions. Section 3 furnishes 12 examples. Section 4 describes a specific algorithm. Section 5 applies the framework to machine learning problems. Section 6 discusses rates of convergence.\\

\state Terminology. A ball under norm $\|\cdot\|$ is denoted by $\ball(\bar x,\rho) = \{x\in\reals^n~|~\|x-\bar x\|\leq \rho\}$. A function $f:\reals^n\to \Reals$ has a {\em domain} $\dom f = \{x\in \reals^n~|~f(x)<\infty\}$ and an {\em epigraph} $\epi f = \{(x,\alpha)\in \reals^{n+1}~|~f(x)\leq \alpha\}$. The function is {\em lower semicontinuous} (lsc) if $\epi f$ is closed and is {\em convex} if $\epi f$ is convex. It is {\em proper} if $\epi f$ is nonempty and $f(x)>-\infty$ for all $x\in \reals^n$.
It is {\em continuous relative to} $C\subset\reals^n$ if $x^\nu\in C\to x\in C$ implies that $f(x^\nu)\to f(x)$. It is {\it locally Lipschitz continuous} (lLc) at $\bar x$ when there are $\delta \in (0,\infty)$ and $\kappa \in [0,\infty)$ such that $|f(x) - f(x')| \leq \kappa\|x-x'\|_2$ whenever $x,x'\in \ball(\bar x,\delta)$. If $f$ is lLc at every $\bar x\in\reals^n$, then $f$ is lLc. A mapping $F:\reals^n\to \reals^m$ is lLc (at $\bar x$) if its component functions $f_1, \dots, f_m$ are lLc (at $\bar x$). Functions and mappings are {\em smooth} (at $\bar x$) if they are continuously differentiable (at $\bar x$). Moreover, $\inf f = \inf\{f(x)~|~x\in\reals^n\}$ and $\epsilon\mbox{-}\nargmin f = \{x\in \dom f~|~f(x) \leq \inf f+\epsilon\}$, which is simply written as $\nargmin f$ if $\epsilon=0$. A {\em lower level-set} $\{f \leq \alpha\} = \{x\in\reals^n~|~f(x) \leq \alpha\}$. The {\em convex hull} of a set $C$ is denoted by $\con C$. The $j$th component of a vector $x\in\reals^n$ is typically indicated by $x_j$. We let $\bdry C$ be the {\em boundary}\footnote{We exclude points in the closure of $C$ that are not in $C$, which is most natural in the present context.} of a set $C\subset\reals^n$, i.e., $\bdry C = C\setminus \nt C$, where $\nt C$ is the interior of $C$. We adopt the usual rules for extended arithmetic, including $\infty - \infty = \infty$.

Sequences of points, sets, and so forth are usually indexed by superscript $\nu\in \nats = \{1, 2, \dots \}$. The collection of subsequences of $\nats$ is denoted by  $\cN_\infty^\grill$, with convergence of $\{x^\nu, \nu\in\nats\}$ to $x$ along a subsequence $N\in \cN_\infty^\grill$ being denoted by $x^\nu \Nto x$. The symbols $\exists$ and $\forall$ mean ``there exist'' and ``for all,'' respectively. The {\em inner limit} of a sequence of sets $\{C^\nu\subset\reals^n, \nu\in\nats\}$ is $\nInnLim C^\nu = \{x\in \reals^n~|~\exists x^\nu \in C^\nu \to x\}$. The {\em outer limit} is $\nOutLim C^\nu = \{x\in \reals^n~|~\exists N\in \cN_\infty^\grill \mbox{ and } x^\nu \in C^\nu \Nto x\}$. Thus, $C^\nu$ {\em set-converges} to $C\subset\reals^n$, denoted by $C^\nu \sto C$, if $\nOutLim C^\nu \subset C \subset \nInnLim C^\nu$. The functions $\{f^\nu:\reals^n\to \Reals, \nu\in\nats\}$ {\em epi-converge} to $f:\reals^n\to \Reals$, denoted by $f^\nu \eto f$, if $\epi f^\nu \sto \epi f$. This takes place if and only if
\begin{align}
  &\forall x^\nu\to x, ~\nliminf f^\nu(x^\nu) \geq f(x)\label{eqn:liminfcondition}\\
  &\forall x, ~\exists x^\nu\to x \mbox{ with } \nlimsup f^\nu(x^\nu) \leq f(x).\label{eqn:limsupcondition}
\end{align}
A set-valued mapping $T:\reals^n\tto \reals^m$ has subsets of $\reals^m$ as its ``values'' and a {\em graph} written as $\gph T = \{(x,y)\in \reals^{n+m}~|~y\in T(x)\}$. A sequence of set-valued mappings $\{T^\nu:\reals^n\tto \reals^m, \nu\in\nats\}$ {\em converges graphically} to $T:\reals^n\tto\reals^m$, denoted by $T^\nu \gto T$, if $\gph T^\nu \sto \gph T$. The set of solutions to the generalized equation $v \in T(x)$ is written as $T^{-1}(v) = \{x \in \reals^n ~|~ v \in T(x)\}$.

We denote the normal cone to $C\subset\reals^n$ at $x\in C$ by $N_C(x)$ and let $N_C(x) = \emptyset$ for $x\not\in C$. Likewise, the set of subgradients of $f:\reals^n\to \Reals$ at a point $x$ with $f(x)\in\reals$ is $\partial f(x)$; $\partial f(x) = \emptyset$ when $f(x)\not\in\reals$. These quantities are understood in the general (Mordukhovich) sense; see \cite[Chap. 6 and 8]{VaAn}.

\section{Consistent Approximations}

Parallel to the actual problem \eqref{eqn:actualpr}, we define the {\em approximating problems}
\begin{equation}\label{eqn:approxpr}
\bigg\{  \nnmin_{x\in \reals^n} ~\phi^\nu(x) = \iota_{X^\nu}(x) + h^\nu\big(F^\nu(x)\big), ~~~~\nu\in\nats\bigg\}
\end{equation}
where $X^\nu\subset\reals^n$, $h^\nu:\reals^m\to \Reals$, and $F^\nu:\reals^n\to \reals^m$ are approximations of the corresponding components of the actual problem. We seek to confirm that small discrepancies between the respective components of the approximating and actual problems indeed translate into small errors in the solutions of the approximating problems relative to those of the actual problem. It is clear that pointwise convergence $\phi^\nu(x)\to \phi(x)$ for all $x\in \reals^n$ would not suffice; see \cite[Figure 7-1]{VaAn}. Uniform convergence of $\phi^\nu$ to $\phi$ on $\reals^n$ ensures that cluster points of $\{x^\nu \in \nargmin \phi^\nu, \nu\in\nats\}$ indeed are minimizers of $\phi$, but such a requirement is too stringent, especially when $X^\nu$ is different than $X$ or $\dom h^\nu$ is different than $\dom h$. We also would like to address convergence of stationary points for the approximating problems to those of the actual problem and then neither pointwise nor uniform convergence of $\phi^\nu$ to $\phi$ is sufficient.

We make the following {\em universal} assumptions for the remainder of the paper:
\[
X,X^\nu \mbox{ are nonempty and closed}; ~~~h,h^\nu \mbox{ are proper, lsc, and convex}; ~~~F,F^\nu \mbox{ are lLc.}
\]
This limits the scope to a well-structured class of problems for which convenient calculus rules can be brought in, while still addressing a vast array of applications. In particular, \cite[Thm. 10.46, Exer. 10.52]{VaAn} lead to optimality conditions for the actual and approximating problems.

If both $h$ and $F$ were smooth, then $0 \in \nabla F(x)^\top \nabla h(F(x)) + N_X(x)$ is a necessary condition for a local minimizer of $\phi$ by \cite[Thm. 6.12]{VaAn}. Here, $\nabla F(x)$ is the Jacobian matrix of $F$ at $x$ with rows corresponding to the gradients of $f_1, \dots, f_m$ at $x$. This condition is equivalently written as
\[
0 = F(x) - z, ~~~~ 0 = \nabla h(z) - y, ~~~~ 0 \in \nsum_{i=1}^m y_i \nabla f_i(x) + N_X(x),
\]
where $y = (y_1, \dots, y_m)\in\reals^m$ and $z\in\reals^m$ are {\em multiplier vectors}. This condition naturally extends to nonsmooth $h$ and $F$ by essentially ``replacing'' gradients by subgradients. Thus, we express optimality conditions for the actual and approximating problems concisely as the generalized equations
\[
0 \in S(x,y,z) ~~~~~~\mbox{ and } ~~~~~~ 0 \in S^\nu(x,y,z)
\]
using the set-valued mappings $S,S^\nu:\reals^{n+2m}\tto \reals^{2m+n}$ given by
\begin{align*}
S(x,y,z) &  = \big\{F(x) - z\big\}  \times \big\{ \partial h(z)- y\big\} \times  \Big( \nsum_{i=1}^m y_i\con \partial f_i(x) + N_X(x)\Big)\\
S^\nu(x,y,z) & = \big\{F^\nu(x) - z\big\} \times \big\{ \partial h^\nu(z)- y\big\} \times \Big( \nsum_{i=1}^m y_i\con \partial f^\nu_i(x) + N_{X^\nu}(x)\Big).
\end{align*}
If $x$ satisfies $0\in S(x,y,z)$ for some $y\in\reals^m$ and $z\in \reals^m$, then $x$ is a {\em stationary point} for the actual problem, with similar terminology being adopted for the approximating problems.
The set of subgradients $\partial f_i(x)$  is convex  when $f_i$ is epi-regular\footnote{A function $f$ is epi-regular at $x$ if its epigraph is Clarke regular \cite[Def. 6.4]{VaAn} at $(x,f(x))$.} at $x$, which, for example, is the case if $f_i$ is smooth at $x$ or if $f_i$ is convex \cite[Thm. 8.30]{VaAn}. Under such circumstances, ``con'' in the expression for $S$ is superfluous.  Since $f_i$ is lLc, $\con \partial f_i(x)$ coincides with the set of Clarke subgradients \cite[Thm. 8.49]{VaAn} of $f_i$ at $x$. Aligned with the approach in \cite{VaAn}, we shy away from this terminology. We view the convexification in the optimality condition as stemming from a slight relaxation of otherwise valid optimality conditions and resulting in an additional ``buffer'' to absorb approximations.

The optimality conditions generalize many familiar ones. For example, $0 \in S(x,y,z)$ if and only if the KKT condition hold for the actual problem with $X = \reals^n$, $h(z) = z_1 + \sum_{i=2}^s \iota_{\{0\}}(z_i) + \sum_{i=s+1}^m \iota_{(-\infty, 0]}(z_i)$, and $F$ being smooth. We see this by deriving $N_X(x) = \{0\}$, $\partial h(z) = \{1\} \times N_{\{0\}^{s-1}}(z_2, \dots, z_s) \times N_{(-\infty, 0]^{m-s}}(z_{s+1}, \dots, z_m)$, and
\[
\nsum_{i=1}^m y_i \con \partial f_i(x) + N_X(x) = \nsum_{i=1}^m y_i \nabla f_i(x)
\]
so that $y_1 = 1$, $y_2, \dots, y_s$ are unrestricted in sign, and, for $i=s+1, \dots, m$, $y_i=0$ if $f_i(x) <0$ and $y_i \geq 0$ if $f_i(x) = 0$. Since even the KKT condition requires a constraint qualification (such as the Mangasarian-Fromovitz), it is not surprising that $0 \in S(x,y,z)$ likewise requires a qualification for it to be a necessary optimality condition. This is formalized in the next proposition.

\begin{proposition}{\rm (optimality condition).}\label{prop:optimality}
Suppose that the following  qualification holds at $x^\star$:
\begin{equation}\label{eqn:qualGenOptCondComposite}
y\in N_{\dom h}\big(F(x^\star)\big) ~~ \mbox{ and }~ ~  0 \in \nsum_{i=1}^m y_i\con \partial f_i(x^\star) + N_X(x^\star)~~~\Longrightarrow~~~ y=0.
\end{equation}
If $x^\star$ is a local minimizer of \eqref{eqn:actualpr}, then 
\[
0\in S(x^\star,y^\star,z^\star)
\]
for some $y^\star\in\reals^m$ and $z^\star\in \reals^m$, with this optimality condition being equivalent to $0 \in \partial \phi(x^\star)$ under the additional assumptions that $X$ is Clarke regular at $x^\star$ and, for each $y\in \partial h(F(x^\star))$ and $i = 1, \dots, m$, $y_i f_i$ is epi-regular at $x^\star$.
\end{proposition}
\state Proof. Guided by \cite[Exer. 10.52]{VaAn}, the first conclusion follows straightforwardly when recalling that $l_F(x,y) = \langle F(x), y\rangle$ implies
\[
\partial_x l_F(x, y) \subset \nsum_{i=1}^m \partial (y_i f_i)(x) \subset \nsum_{i=1}^m \con \big\{\partial (y_i f_i)(x)\big\} = \nsum_{i=1}^m y_i \con \partial f_i(x)
\]
via \cite[Cor. 10.9, Thm. 9.61]{VaAn}.

The second conclusion about equivalence with $0\in \partial \phi(x^\star)$ holds by the following argument. For $\tilde h:\reals^n\times\reals^m\to \Reals$ and $\tilde F:\reals^n\to \reals^{n+m}$ given by
\[
\tilde h(z_0,z) = \iota_X(z_0) + h(z)~~~\mbox{ and } ~~~ \tilde F(x) = \big(x, F(x)\big), 
\]
we obtain that $\phi(x) = \tilde h(\tilde F(x))$. Clearly, $\tilde F$ is lLc and $\tilde h$ is lsc and proper. Let $l_{\tilde F}(x,(w,y)) = \langle \tilde F(x), (w,y)\rangle$. Under a qualification, \cite[Thm. 6.23]{primer} confirms that
\[
\partial \phi(x^\star) = \bigcup_{(w,y)\in \partial \tilde h(\tilde F(x^\star))} \partial_x l_{\tilde F}\big(x^\star,(w,y)\big)
\]
as long as $\tilde h$ is epi-regular at $\tilde F(x^\star)$ and $l_{\tilde F}(\cdot\,,(w,y))$ is epi-regular at $x^\star$ for all $(w,y) \in \partial \tilde h(\tilde F(x^\star))$. A closer examination shows that \eqref{eqn:qualGenOptCondComposite} ensure that the qualification of the theorem indeed holds. By \cite[Prop. 4.63]{primer}, we conclude that $\tilde h$ is epi-regular at $\tilde F(x^\star)$ because $X$ is Clarke regular at $x^\star$ and $\partial h(F(x^\star)) \neq \emptyset$. Moreover, $l_{\tilde F}(\cdot\,,(w,y))$ is epi-regular at $x^\star$ because $y_i f_i$ is assumed to be epi-regular at $x^\star$; see \cite[Exam. 4.70]{primer}. Thus, it only remains to untangle the expression for $\partial \phi(x^\star)$. We obtain that
\[
0 \in \partial \phi(x^\star) ~~~\Longleftrightarrow ~~~ \exists y \in \partial h\big(F(x^\star)\big), ~~ 0 \in \nsum_{i=1}^m  \partial (y_i f_i)(x^\star) + N_X(x^\star),
\]
where we again appeal to \cite[Exam. 4.70]{primer} to compute $\partial_x l_{\tilde F}(x^\star,(w,y))$. Since $y_if_i$ is epi-regular at $x^\star$, $\partial (y_if_i)(x^\star)$ is a convex set and we conclude that
\[
\partial (y_i f_i)(x^\star) = \con \big\{\partial (y_i f_i)(x^\star)\big\} = y_i \con \partial f_i(x^\star),
\]
which implies the assertion.\eop

In the absence of regularity of the kind invoked in the second part of the proposition, the optimality condition $0 \in \partial \phi(x^\star)$ could be strictly stronger than $0 \in S(x^\star, y^\star, z^\star)$. For example, let $X = \reals$, $h(z) = z_1 + z_2$, and $F(x) = (\max\{-x, x/2\}, \min\{0,-x\})$. Then, $\phi(x) = \iota_X(x) + h(F(x)) = \max\{-x,-x/2\}$ with $0 \not\in\partial \phi(0) = [-1,-1/2]$. In contrast, $0\in S(0,y^\star,z^\star)$ for $y^\star = (1,1)$ and $z^\star = (0,0)$ because 
\[
S(0,y^\star,z^\star) = \big\{F(0) - z^\star\big\}  \times \big\{ \partial h(z^\star) - y^\star\big\} \times  \big( y_1^\star \con [-1, 1/2] + y_2^\star \con \{-1, 0\} \big).
\]
Generally, the advantage of the optimality condition $0 \in S(x,y,z)$ is its explicit form in terms of the ``primitives'' of the actual problem \eqref{eqn:actualpr}. It is therefore much more computationally accessible than $0\in \partial \phi(x)$, which can be checked numerically only in special cases.

Proposition \ref{prop:optimality} likewise confirms that $0\in S^\nu(x,y,z)$ is a necessary optimality condition for an approximating problem under a qualification parallel to \eqref{eqn:qualGenOptCondComposite}. We are now in a position to state what it means for the approximating problems, paired with their optimality conditions, to be well-justified surrogates of the actual problem and its optimality condition.

\begin{definition}{\rm (consistent approximations).}
The pairs $\{(\phi^\nu,S^\nu), \nu\in \nats\}$ are {\em consistent approximations} of $(\phi,S)$ when
\[
\phi^\nu \eto \phi~~~~\mbox{ and }~~~~S^\nu\gto S.
\]
If the graphical convergence is relaxed to merely $\nOutLim (\gph S^\nu) \subset \gph S$, then $\{(\phi^\nu,S^\nu), \nu\in \nats\}$ are {\em weakly consistent approximations} of $(\phi,S)$.
\end{definition}

Under consistency, even in its weak form, we are on solid ground because the approximating problems indeed produce approximating solutions of the actual problem as formalized next using the following notation. For $T:\reals^n\tto \reals^m$, the {\it set of $\delta$-solutions} to the generalized equation $0 \in T(x)$ is defined as
\[
T^{-1}\big( \ball(0, \delta)  \big) = \bigcup_{v\in \ball(0, \delta)} T^{-1}(v) = \bigcup_{v\in \ball(0, \delta)} \big\{x\in\reals^n~\big|~v\in T(x)\big\}.
\]
Thus, a $\delta$-solution $\hat x\in T^{-1}( \ball(0, \delta)  )$ ``almost'' satisfies $0 \in T(x)$ in the sense that $v \in T(\hat x)$ for some $v$ with $\|v\| \leq \delta$. We observe that such near-solutions may depend on the choice of norm $\|\cdot\|$.

\begin{proposition}\label{tconvConAlgo}{\rm (consequences of consistency).} Suppose that $\{(\phi^\nu,S^\nu), \nu\in \nats\}$ are weakly consistent approximations of $(\phi,S)$, the tolerances $\epsilon^\nu$ and $\delta^\nu$ vanish, and $\alpha<\beta$. Then, the following hold:
\begin{enumerate}[{\rm (a)}]
\item $\nOutLim \big(\epsilon^\nu\mbox{-}\nargmin \phi^\nu\big) \subset \nargmin \phi$, provided that $\dom \phi \neq \emptyset$.

\item $\nOutLim \{\phi^\nu \leq \alpha\} \subset \{ \phi \leq \alpha \} \subset \nInnLim \{\phi^\nu \leq \beta\}$.

\item $\nOutLim (S^\nu)^{-1}\big(\ball(0,\delta^\nu)\big) \subset S^{-1}(0)$.

\end{enumerate}
Moreover, if the approximations are consistent (and not merely weakly consistent), then there are $\gamma^\nu \to 0$ such that $\nInnLim (S^\nu)^{-1}(\ball(0,\gamma^\nu)) \supset S^{-1}(0)$.

\end{proposition}
\state Proof. This is a compilation of well-known facts; see \cite[Thm. 7.31]{VaAn} for (a); \cite[Prop. 7.7]{VaAn} for (b); and \cite[Thm. 5.37]{VaAn} for (c) and the final statement.\eop

Consistency indeed guarantees a comprehensive sense of approximation tailored to the need of minimization problems. From (a) in the proposition, we see that every cluster point of a sequence of near-minimizers of the approximating problems is a minimizer of the actual problem provided that the tolerances vanish. That is, if $x^\nu \in \epsilon^\nu\mbox{-}\nargmin \phi^\nu \Nto \bar x$ along a subsequence $N\in \cN_\infty^\grill$, then $\bar x\in \nargmin \phi$. This is a valuable property, but remains somewhat conceptual in nonconvex settings because there might be no practical algorithm for computing near-minimizers of the approximating problems. Item (b) shows that the lower level-sets of $\phi^\nu$ eventually become essentially indistinguishable from those of $\phi$. For example, if $\phi^\nu(x^\nu) \leq \alpha$ and $x^\nu\to \bar x$, then $\phi(\bar x) \leq \alpha$ by the first inclusion in (b). The second inclusion in (b) establishes that any $\bar x$ with $\phi(\bar x) \leq \alpha$ can be approached by points $x^\nu$ with $\phi^\nu(x^\nu) \leq \beta$. Thus, $\phi^\nu$ accurately ``predicts'' the values of $\phi$ as $\nu\to \infty$.

Item (c) confirms that any cluster point $(\bar x, \bar y, \bar z)$ of a sequence $\{(x^\nu,y^\nu,z^\nu),\, \nu\in \nats\}$, with $0 \in S^\nu(x^\nu,y^\nu,z^\nu)$, satisfies $0 \in S(\bar x,\bar y,\bar z)$. This holds even if the points $(x^\nu,y^\nu,z^\nu)$ are computed with tolerances $\delta^\nu$ as long as they vanish. The approximating problems and their optimality conditions are often constructed such that computing nearly stationary points is indeed possible in finite time; see Section 3 for examples. Combining the facts from (b) and (c), we summarize: If $N\in \cN_\infty^\grill$ and $(x^\nu,y^\nu,z^\nu)\in (S^\nu)^{-1}\big(\ball(0,\delta^\nu)\big) \Nto (\bar x, \bar y, \bar z)$, then
    \[
    (\bar x, \bar y, \bar z) \in S^{-1}(0) ~~\mbox{ and }~~ \phi(\bar x) \leq \nliminf_{\nu\in N} \phi^\nu(x^\nu).
    \]
Thus, $\bar x$ is stationary for the actual problem, with an objective function value at least as good as predicted by the approximations.

These consequences occur under {\em weak} consistency. Passing to consistency, we achieve enhanced robustness in the sense that every stationary point of the actual problem can be approached by nearly stationary points of the approximating problems. In fact, the proposition guarantees that for some tolerances $\gamma^\nu\to 0$, one has $(S^\nu)^{-1}(\ball(0,\gamma^\nu)) \sto S^{-1}(0)$.

The proposition justifies the following broad algorithmic framework for solving the actual problem.

\medskip

\state Consistent Approximation Algorithm.

\begin{description}

  \item[Data.] ~\,~$\delta^\nu \geq 0$, with $\delta^\nu\to 0$.

  \item[Step 0.]  Set $\nu = 1$.

  \item[Step 1.]  Minimize $h^\nu \comp F^\nu$ over $X^\nu$ until one obtains $x^\nu$, with corresponding $(y^\nu,z^\nu)$, satisfying
  \[
  (x^\nu,y^\nu,z^\nu)\in (S^\nu)^{-1}\big(\ball(0,\delta^\nu)\big).
\]

\item[Step 2.] Replace $\nu$ by $\nu +1$ and go to Step 1.
\end{description}

\medskip

A multitude of implementation details remain unsettled. What approximations are most suitable for the actual problem at hand? What subroutine should be used in Step 1 and with what $\delta^\nu$? We list many examples of approximations in Section 3 and describe a concrete algorithm in Section 4.

Verification of consistency is supported by sufficient conditions for epi-convergence of $\phi^\nu$ to $\phi$ and graphical convergence of $S^\nu$ to $S$. There are several known results. For example, if $X = X^\nu = \reals^n$, $F = F^\nu$ are smooth with $\nabla F$ having rank $m$,  and $h^\nu\eto h$, then $\phi^\nu\eto \phi$ by \cite[Thm. 3.2]{BurkeHoheisel.13}; see also \cite[Chap. 5 and 7]{VaAn}. The following conditions appear especially versatile by allowing for approximations of $X$, $h$, and $F$ under relatively mild assumptions.

\begin{theorem}{\rm (sufficient conditions for consistent approximations).}\label{tConstentApproxComp}
Suppose that $h^\nu \eto h$, $X^\nu \sto X$ with these sets being convex, and, for $i=1, \dots, m$, one has the property:
\begin{equation}\label{eqn:suffFconv}
  \begin{rcases}
    x^\nu \in X^\nu \to x\in X\\
    v^\nu \in \partial f_i^\nu(x^\nu)
  \end{rcases}
  ~\Longrightarrow ~
  \begin{cases}
    f_i^\nu(x^\nu) \to f_i(x)\\
    \{v^\nu, \nu\in\nats\} \mbox{ is bounded with all its cluster points in } \con \partial f_i(x).
  \end{cases}
\end{equation}
Then, each one of the following conditions is sufficient for $\{(\phi^\nu,S^\nu), \nu\in\nats\}$ to be weakly consistent approximations of $(\phi,S)$:
\begin{enumerate}[{\rm (a)}]

  \item {\rm (real-valuedness):} $h$ is real-valued.

  \item {\rm (pointwise convergence):} $h^\nu(z)\to h(z)$ for all $z\in \bdry (\dom h)$ and, in addition, $F^\nu(x) = F(x)$ for all $x\in X$ and $X \subset X^\nu$ for each $\nu$.

  \item {\rm (monotonicity):} For each $\nu$: $h^\nu(z)\leq h^\nu(\bar z)$ when $z \leq \bar z$, $h^\nu(z) \leq h(z)$ for all $z$, $X \subset X^\nu$, and $F^\nu(x) \leq F(x)$ for all $x\in X$.

  \item {\rm (interior points):} $h$ is continuous relative to $\dom h$, $X \subset X^\nu$ for each $\nu$, and for all $x\in X$ with $F(x)\in \bdry(\dom h)$, there is $x^\nu\in X \to x$ such that $F(x^\nu) \in \nt(\dom h)$.

  \item {\rm (separability):} For $\sum_{k=1}^r m_k = m$, lLc $F_k^\nu, F_k:\reals^n\to \reals^{m_k}$, and proper, lsc, and convex $h^\nu_k, h_k:\reals^{m_k}\to \Reals$, $k=1, \dots, r$, one can express $F^\nu(x) = (F^\nu_1(x), \dots, F^\nu_r(x))$, $F(x) = (F_1(x), \dots, F_r(x))$, $h^\nu(z) = \nsum_{k=1}^r h_k^\nu(z_k)$, and $h(z) = \nsum_{k=1}^r h_k(z_k)$ with $h_k$ being either real-valued or $h_k$ satisfying $h^\nu_k(z_k) \to h_k(z_k)$ for all $z_k\in \bdry (\dom h_k)$ and $F^\nu_k(x) = F_k(x)$ for all $x\in X$. Moreover, $X \subset X^\nu$.

\end{enumerate}
We also have the refinements:

\noindent If $X^\nu = X$, then each of (a)-(e) remains sufficient without the convexity assumption on $X$.

\noindent If $f_i^\nu = f_i$, then \eqref{eqn:suffFconv} holds automatically.

\noindent If $f_i^\nu, f_i$ are convex and $f_i^\nu\eto f_i$, then \eqref{eqn:suffFconv} holds automatically.

\noindent If $X^\nu = X$ and (c) holds, then $f_i^\nu(x^\nu) \to f_i(x)$ in  \eqref{eqn:suffFconv} is relaxed to $\nliminf f_i^\nu(x^\nu) \geq f_i(x)$. 

\noindent If $X^\nu = X$ and $f_i^\nu = f_i$, then each of (a)-(e) is sufficient for consistent approximations.

\noindent If $f_i^\nu, f_i$ are smooth, then \eqref{eqn:suffFconv} is equivalent to $f_i^\nu(x^\nu) \to f_i(x)$ and $\nabla f_i^\nu(x^\nu) \to \nabla f_i(x)$ as $x^\nu \in X^\nu \to x$.

\noindent If $f_i^\nu, f_i$ are smooth, then each of (a)-(e) is sufficient for consistent approximations.
\end{theorem}
\state Proof. We start by establishing $\phi^\nu\eto \phi$. To prove the liminf-condition \eqref{eqn:liminfcondition} for $\phi^\nu$ and $\phi$, let $x^\nu\to x$. Since neither $\iota_X(x)$ nor $h(F(x))$ is equal to $-\infty$,
\[
\nliminf \Big( \iota_{X^\nu}(x^\nu) + h^\nu\big(F^\nu(x^\nu)\big)\Big) \geq \nliminf \iota_{X^\nu}(x^\nu) + \nliminf h^\nu\big(F^\nu(x^\nu)\big).
\]
Certainly, $X^\nu\sto X$ implies that $\iota_{X^\nu}\eto \iota_X$ and then also $\nliminf \iota_{X^\nu}(x^\nu) \geq \iota_X(x)$. If $x^\nu\not\in X^\nu$, then $\iota_{X^\nu}(x^\nu)=\infty$. Thus, we assume without loss of generality that $x^\nu\in X^\nu$. This implies that $x\in X$ because $X^\nu\sto X$. By assumption, $F^\nu(x^\nu) \to F(x)$ and $h^\nu \eto h$. These facts ensure that $\nliminf h^\nu(F^\nu(x^\nu))$ $\geq$ $h(F(x))$. When combined with the earlier inequalities, this relation establishes the liminf-condition \eqref{eqn:liminfcondition} for $\phi^\nu$ and $\phi$, which thus holds without any of the additional conditions (a)-(e).

Next, we establish the limsup-condition \eqref{eqn:limsupcondition} for $\phi^\nu$ and $\phi$. Let $x\in\reals^n$ be arbitrary. If $x\not\in X$ or $F(x)\not\in \dom h$, then $\phi(x) = \infty$ and the limsup-condition holds trivially. If $x\in X$ and $F(x)\in \dom h$, then we argue as follows:

First, suppose that (a) holds. Since $X^\nu\sto X$, there exists $x^\nu\in X^\nu \to x$. Moreover, $h^\nu$ converges uniformly to $h$ on compact sets because $h$ is convex and real-valued; see \cite[Thm. 7.17]{VaAn}. Consequently, $h^\nu(F^\nu(x^\nu))\to h(F(x))$ and one has
\[
\nlimsup \Big( \iota_{X^\nu}(x^\nu) + h^\nu\big(F^\nu(x^\nu)\big)\Big) \leq \nlimsup \iota_{X^\nu}(x^\nu) + \nlimsup h^\nu\big(F^\nu(x^\nu)\big) = h\big(F(x)\big).
\]

Second, suppose that (c) holds. Since $X\subset X^\nu$, we construct $\{x^\nu= x, \nu\in\nats\}$ and simply need to prove that $\nlimsup h^\nu(F^\nu(x))\leq h(F(x))$. This holds because $F^\nu(x) \leq F(x)$ and then
$h^\nu(F^\nu(x)) \leq h^\nu(F(x)) \leq h(F(x))$ for all $\nu\in\nats$.

Third, suppose that (d) holds. We consider two cases. Suppose that $F(x) \in \nt (\dom h)$. Then, construct $\{x^\nu = x, \nu\in \nats\}$ and note that $h^\nu(F^\nu(x)) \to h(F(x))$ because $h^\nu$ converges uniformly to $h$ on bounded subsets of the interior of $\dom h$; see \cite[Thm. 7.17]{VaAn}. Next, suppose that $F(x) \in \bdry (\dom h)$. Then, by condition (d), there are $\bar x^k\in X\to x$ with $F(\bar x^k) \in \nt (\dom h)$. Consequently, there exist $\delta_k \in (0,\infty) \to 0$ such that $\ball(F(\bar x^k), \delta_k) \subset \nt (\dom h)$ for all $k$. Let $\epsilon_k \in (0,\infty) \to 0$. Fix $k$. Since $\bar x^k\in X \cap X^\nu$, $F^\nu(\bar x^k) \to F(\bar x^k)$ as $\nu\to \infty$. Let $\nu_k \in \nats$, with $\nu_k\geq k$, be such that $F^\nu(\bar x^k) \in \ball(F(\bar x^k),\delta_k)$ for all $\nu\geq \nu_k$. Since $h^\nu$ converges uniformly to $h$ on $\ball(F(\bar x^k),\delta_k)$ by virtue of that set being in the interior of $\dom h$ (again cf. \cite[Thm. 7.17]{VaAn}), there is $\bar \nu_k\geq \nu_k$ such that $h^\nu(F^\nu(\bar x^k)) \leq h(F(\bar x^k)) + \epsilon_k$ for all $\nu\geq \bar \nu_k$. We repeat these arguments for all $k\in \nats$. Now, we construct $x^\nu = \bar x^1$ for $\nu = 1,\dots, \bar \nu_2$ and  $x^\nu = \bar x^k$ for $\bar\nu_k< \nu \leq \bar \nu_{k+1}$. Then, $x^\nu \to x$ and $h^\nu(F^\nu(x^\nu)) \leq h(F(x^\nu)) + \epsilon_k$ for all $\bar\nu_k< \nu \leq \bar \nu_{k+1}$. As $\nu\to \infty$, $k$ tends to $\infty$ as well in this expression. Thus, its right-hand side converges to $h(F(x))$; recall that $h$ is continuous on its domain. We have shown that $\nlimsup h^\nu(F^\nu(x^\nu)) \leq h(F(x))$.

Fourth, suppose that (e) holds. Since $X\subset X^\nu$, one has $x\in X^\nu$ for all $\nu$ and we construct $\{x^\nu= x, \nu\in\nats\}$. Certainly, $F^\nu(x) \to F(x)$ and
\[
\nlimsup \Big( \iota_{X^\nu}(x^\nu) + \nsum_{k=1}^r h_k^\nu\big(F_k^\nu(x^\nu)\big)\Big) \leq \nlimsup \iota_{X^\nu}(x) + \nsum_{k=1}^r \nlimsup h_k^\nu\big(F^\nu_k(x)\big).
\]
Suppose that $h_k$ is real-valued. Then, $h_k^\nu$ converges uniformly to $h_k$ on compact sets because $h_k$ is convex; see \cite[Thm. 7.17]{VaAn}. This implies that $h_k^\nu(F^\nu_k(x)) \to h_k(F_k(x))$. Alternatively, suppose that  $h^\nu_k(z_k) \to h_k(z_k)$ for all $z_k\in \bdry (\dom h_k)$ and $F^\nu_k(x') = F_k(x')$ for all $x'\in X$. Since $F(x) \in \dom h$, we have $F_k(x) \in \dom h_k$. If $F_k(x)\in \bdry (\dom h_k)$, then, by assumption, $h^\nu_k(F_k(x)) \to h_k(F_k(x))$. If $F_k(x)\in \nt (\dom h_k)$, then $h^\nu_k(F_k(x)) \to h_k(F_k(x))$ as well because these functions are convex. In summary,
\[
\nsum_{k=1}^r \nlimsup h_k^\nu\big(F^\nu_k(x)\big) = \nsum_{k=1}^r h_k\big(F_k(x)\big)
\]
and the needed limsup-condition holds.

Fifth, suppose that (b) holds. This is a special case of (e) with $r = 1$.

We have confirmed that each one of (a)-(e) suffices for $\phi^\nu\eto \phi$ and, in fact, this holds even without the convexity of $X^\nu$ and $X$.

We next turn to the relation between $\gph S^\nu$ and $\gph S$ and start by showing $\nOutLim (\gph S^\nu) \subset \gph S$. Let $(\bar x, \bar y, \bar z, \bar u, \bar v, \bar w)\in \nOutLim (\gph S^\nu)$. Then, there are $N\in \cN_\infty^\grill$, $x^\nu\Nto \bar x$, $y^\nu\Nto \bar y$, $z^\nu\Nto \bar z$, $u^\nu\Nto \bar u$, $v^\nu\Nto \bar v$, and $w^\nu\Nto \bar w$ with
$(x^\nu, y^\nu, z^\nu, u^\nu, v^\nu, w^\nu) \in \gph S^\nu$. Consequently,
\[
u^\nu = F^\nu(x^\nu) - z^\nu, ~~~ v^\nu \in \partial h^\nu(z^\nu)- y^\nu, ~~~ w^\nu \in \nsum_{i=1}^m y_i^\nu \con \partial f_i^\nu(x^\nu) + N_{X^\nu}(x^\nu).
\]
The existence of $w^\nu$ implies that $x^\nu\in X^\nu$ because otherwise $N_{X^\nu}(x^\nu)$ would have been an empty set. Since $X^\nu\sto X$, we also have $\bar x\in X$. Thus,  $F^\nu(x^\nu)\Nto F(\bar x)$ and we conclude that $\bar u = F(\bar x) - \bar z$. Attouch's theorem \cite[Thm. 12.35]{VaAn} states that $h^\nu \eto h$ implies $\partial h^\nu \gto \partial h$. Consequently, $\bar v \in \partial h(\bar z) - \bar y$.

The inclusion for $w^\nu$ implies that there exist
\begin{equation}\label{eqn:constProofw}
a_i^\nu \in \con \partial f_i^\nu(x^\nu), ~i=1, \dots, m,~ ~~\mbox{ such that } ~~~ w^\nu - \nsum_{i=1}^m y_i^\nu a_i^\nu \in N_{X^\nu}(x^\nu).
\end{equation}
By Caratheodory's theorem, there are $\lambda_{ij}^\nu \geq 0$ and $b_{ij}^\nu\in \partial f_i^\nu(x^\nu)$, $j=0, 1, \dots, n$, such that $\sum_{j=0}^n \lambda_{ij}^\nu = 1$ and $a_i^\nu = \sum_{j=0}^n \lambda_{ij}^\nu b_{ij}^\nu$. The sequence $\{(\lambda_{i0}^\nu, \dots, \lambda_{in}^\nu), \nu\in N\}$ is contained in a compact set and thus has a convergent subsequence. By assumption \eqref{eqn:suffFconv}, $\{b_{ij}^\nu, \nu \in N\}$ is bounded with all its cluster points in $\con \partial f_i(\bar x)$. Consequently, there exist $\bar \lambda_i = (\bar \lambda_{i0}, \dots, \bar \lambda_{in})\in \reals^{1+n}$, $\bar b_{i0}, \dots, \bar b_{in}$, and a subsequence of $N$, which we also denote by $N$, such that
\[
(\lambda_{i0}^\nu, \dots, \lambda_{in}^\nu) \Nto \bar\lambda_i \geq 0, ~\mbox{ with } ~\nsum_{j=0}^n \bar \lambda_{ij} = 1,  ~~~\mbox{ and } ~~~ b_{ij}^\nu \Nto \bar b_{ij} \in \con \partial f_i(\bar x).
\]
This means that
\[
a_i^\nu \Nto \bar a_i = \nsum_{j=0}^n \bar \lambda_{ij} \bar b_{ij} \in \con \partial f_i(\bar x).
\]
Via Attouch's theorem  \cite[Thm. 12.35]{VaAn}, $\iota_{X^\nu} \eto \iota_X$ implies $N_{X^\nu} \gto N_X$ because $X^\nu$ and $X$ are convex. Using this fact as well as the observations that $y_i^\nu \Nto \bar y_i$ and $w^\nu - \sum_{i=1}^m y_i^\nu a_i^\nu \Nto \bar w  - \sum_{i=1}^m \bar y_i \bar a_i$, we obtain from \eqref{eqn:constProofw} that $\bar w - \nsum_{i=1}^m \bar y_i \bar a_i \in N_{X}(\bar x)$. Thus, because $\bar a_i \in \con \partial f_i(\bar x)$, one has
\[
\bar w \in \nsum_{i=1}^m \bar y_i \con \partial  f_i(\bar x) + N_{X}(\bar x)
\]
and we conclude that $(\bar u, \bar v, \bar w) \in S(\bar x, \bar y, \bar z)$. We have established that $\nOutLim (\gph S^\nu) \subset \gph S$ and this holds without leveraging any of the conditions (a)-(e).

Under the modified assumption that $X^\nu = X$, but not necessarily convex, we still have $N_{X^\nu} =N_X \gto N_X$ because $N_X$ is outer semicontinuous; see \cite[Prop. 6.6]{VaAn}. Thus, the above argument carries over.

If $f_i^\nu = f_i$, then \eqref{eqn:suffFconv} holds automatically because $f_i$ is lLc and $\partial f_i$ is outer semicontinuous; see \cite[Prop. 8.7]{VaAn}.

If $X^\nu = X$ and (c) holds, then $\nlimsup f_i^\nu(x^\nu) \leq \nlimsup f_i(x^\nu) = f_i(x)$, with the last equality following by continuity of $f_i$.  

If $f_i^\nu,f_i$ are convex and $f_i^\nu\eto f$, then \cite[Prop. 4.18]{primer} establishes that $f_i^\nu(x^\nu) \to f_i(\bar x)$ if $x^\nu \to \bar x$. Cluster points of the sequence $\{v^\nu \in \partial f_i^\nu(x^\nu), \nu\in\nats\}$ must be in $\partial f_i(\bar x)$ by Attouch’s theorem; see \cite[Thm. 12.35]{VaAn}. It only remains to show that the sequence is bounded. We recall that
\[
\partial f(x) = \nargmin_v \big\{ f^*(v) - \langle v,x\rangle \big\}  ~~~ \mbox{ and } ~~~ \partial f^\nu(x) = \nargmin_v \big\{ f^{\nu*}(v) - \langle v,x\rangle \big\}
\]
for every $x\in\reals^n$; see \cite[Prop. 11.3]{VaAn}. Here, $f^*$ and $f^{\nu*}$ are the conjugates of $f$ and $f^\nu$, respectively, which are also proper, lsc, and convex by the Fenchel-Moreau theorem 5.23 in \cite{primer}.  By Wijsman's theorem \cite[Thm. 11.34]{VaAn}, $f^{\nu*} \eto f^*$ and then we also have $g^\nu = f^{\nu*} - \langle \cdot\,,x^\nu\rangle \eto g = f^{*} - \langle \cdot\,,\bar x\rangle$ as seen from \cite[Prop. 4.19(a)]{primer}. Since $\partial f(\bar x)$ is compact by \cite[Prop. 2.54]{primer}, there is $\rho \in (0,\infty)$ such that $\partial f(\bar x) = \nargmin g \subset \ball(0,\rho/2)$. Let $B =  \ball(0,\rho)$. Set $\eta = \inf_{\bdry B} g$. If $\eta<\infty$, then $g(v)<\infty$ for some $v\in \bdry B$ so \cite[Thm. 4.9]{primer} applies, $\nargmin_{\bdry B} g \neq\emptyset$,  and $\eta = g(v^\star)$ for all $v^\star \in \nargmin_{\bdry B} g$. Since $\nargmin g \cap \bdry B = \emptyset$, $g(v^\star)>\inf g$ for $v^\star \in \nargmin_{\bdry B} g$. Thus, regardless of $\eta$ being finite or not, there is $\delta \in (0,\infty)$ such that
\begin{equation}\label{eqn:LpConvMinConvex0}
g(v)\geq \inf g + \delta~~~ \forall v\in \bdry B.
\end{equation}
Since $\partial f(\bar x)$ is nonempty by \cite[Prop. 2.25]{primer}, there is $\bar v \in \nargmin g$. From the definition of epi-convergence, $g^\nu\eto g$ implies that there are points $w^\nu\to \bar v$ such that $\nlimsup g^\nu(w^\nu)\leq g(\bar v)$. Thus, there exists $\nu_0$ such that $g^\nu(w^\nu) \leq \inf g + \third\delta$  and $w^\nu\in B$  for all $\nu\geq \nu_0$. For the sake of contradiction, suppose that there are $N\in\cN_\infty^\grill$ such that $\{v^\nu\not\in B, \nu\in N\}$. Since $v^\nu \in \nargmin g^\nu$, one has $g^\nu(v^\nu) \leq g^\nu(w^\nu)$ for all $\nu\in N$. For $\nu\in N$ with $\nu\geq \nu_0$, let $z^\nu$ be the unique point in $\bdry B$ on the line segment between $v^\nu$ and $w^\nu$, i.e., $z^\nu = (1-\lambda^\nu) v^\nu + \lambda^\nu w^\nu$  for some $\lambda^\nu \in [0,1]$. The convexity inequality (\cite[Prop. 1.10]{primer}) applied to $g^\nu$ implies that
\begin{align}\label{eqn:pConvMinConvex}
g^\nu(z^\nu) & \leq (1-\lambda^\nu) g^\nu(v^\nu) + \lambda^\nu g^\nu(w^\nu)\nonumber\\
&  \leq (1-\lambda^\nu) g^\nu(w^\nu) + \lambda^\nu g^\nu(w^\nu) = g^\nu(w^\nu) \leq \inf g + \third\delta
\end{align}
for $\nu\in N$ with $\nu\geq \nu_0$. Since $\{z^\nu, \nu\in N, \nu\geq \nu_0\}$ is contained in the compact set $\bdry B$, it has a cluster point $\bar z\in \bdry B$. After passing to the corresponding subsequence, which we also denote by $N$, we find that $\nliminf_{\nu\in N} g^\nu(z^\nu) \geq g(\bar z) \geq \inf g + \delta$, where the first inequality follows from the fact that $g^\nu\eto g$ and the second one by \eqref{eqn:LpConvMinConvex0}. Thus, for sufficiently large $\nu\in N$, $g^\nu(z^\nu) \geq \inf g + \twothird\delta$. However, this contradicts \eqref{eqn:pConvMinConvex} and we conclude that $v^\nu\in B$ for all but a finite number of $\nu$.

Under the assumption that $X^\nu = X$ and $f_i^\nu = f_i$, we show that $\nInnLim (\gph S^\nu) \supset \gph S$. Let $(\bar x, \bar y, \bar z, \bar u,
\bar v, \bar w) \in \gph S$. This means that
\[
\bar u = F(\bar x) - \bar z, ~~~ \bar v \in \partial h(\bar z)- \bar y, ~~~ \bar w \in \nsum_{i=1}^m \bar y_i \con \partial f_i(\bar x) + N_{X}(\bar x).
\]
Set $x^\nu = \bar x$, $y^\nu = \bar y$, and $w^\nu = \bar w$. Since $h^\nu\eto h$, it follows by Attouch's theorem
\cite[Thm. 12.35]{VaAn} that $\partial h^\nu \gto \partial h$ and there are $v^\nu \to \bar v$ and $z^\nu\to \bar z$ such that $v^\nu + \bar y \in \partial h^\nu(z^\nu)$. Let $u^\nu = F(\bar x) - z^\nu$. This means that $(x^\nu, y^\nu, z^\nu, u^\nu, v^\nu, w^\nu) \in \gph S^\nu \to (\bar x, \bar y, \bar z, \bar u, \bar v, \bar w)$. Consequently, $(\bar x, \bar y, \bar z, \bar u, \bar v, \bar w) \in \nInnLim (\gph S^\nu)$ and $\nInnLim (\gph S^\nu) \supset \gph S$ holds. In view of the
earlier results, $S^\nu \gto S$ and we conclude that the approximations are consistent.

For the case with smooth $f_i^\nu,f_i$, the claim about \eqref{eqn:suffFconv} is trivial and it only remains to show $\nInnLim (\gph S^\nu) \supset \gph S$. Now, we obtain the simplification
\[
\nsum_{i=1}^m y_i \con \partial f_i(x) = \Big\{\nsum_{i=1}^m y_i \nabla f_i(x)\Big\} = \big\{\nabla F(x)^\top y\big\},
\]
with a similar expression for the approximating functions. Let $(\bar x, \bar y, \bar z, \bar u, \bar v, \bar w) \in \gph S$. This means that $\bar u = F(\bar x) - \bar z$, $\bar v \in \partial h(\bar z)- \bar y$, and $\bar w \in \nabla F(\bar x)^\top \bar y + N_{X}(\bar x)$.
In particular, $\bar x\in X$. Set $y^\nu = \bar y$. Since $h^\nu\eto h$, it follows by Attouch's theorem  \cite[Thm. 12.35]{VaAn} that $\partial h^\nu \gto \partial h$ and there are $v^\nu \to \bar v$ and $z^\nu\to \bar z$ such that $v^\nu + \bar y \in \partial h^\nu(z^\nu)$.

Since $N_{X^\nu} \gto N_X$, as argued above, there's $(x^\nu,t^\nu) \in \gph N_{X^\nu}$ with $x^\nu \in X^\nu \to \bar x$ and $t^\nu \to \bar w - \nabla F(\bar x)^\top \bar y$. Construct $w^\nu = t^\nu +  \nabla F^\nu(x^\nu)^\top \bar y$. We then have $w^\nu \to \bar w$ and $w^\nu - \nabla F^\nu(x^\nu)^\top \bar y \in N_{X^\nu}(x^\nu)$. Also, construct $u^\nu = F^\nu(x^\nu) - z^\nu$, which converges to $\bar u$.
In summary, we have constructed $(x^\nu, y^\nu, z^\nu, u^\nu, v^\nu, w^\nu) \in \gph S^\nu \to (\bar x, \bar y, \bar z, \bar u, \bar v, \bar w)$. This means that $(\bar x, \bar y, \bar z, \bar u, \bar v, \bar w) \in \nInnLim (\gph S^\nu)$ and, thus, $\nInnLim (\gph S^\nu) \supset \gph S$ holds.
\eop

While we believe weak consistency is the natural goal in most applications, in some cases one might be willing to sacrifice some assurances with the benefit of relaxed assumptions.

\begin{corollary}{\rm (approximations without full epi-convergence).}\label{cApprox}
Suppose that $h^\nu \eto h$, $X^\nu \sto X$ with these sets being convex,  and \eqref{eqn:suffFconv} holds  for $i=1, \dots, m$. Then, for $\alpha \in \Reals$ and vanishing $\delta^\nu$, one has
\[
\nOutLim \{\phi^\nu \leq \alpha\} \subset \{ \phi \leq \alpha \}~~~~~\mbox{ and } ~~~~~ \nOutLim (S^\nu)^{-1}\big(\ball(0,\delta^\nu)\big) \subset S^{-1}(0).
\]
If $X^\nu = X$, then the assertions hold without the convexity assumption on $X$.
\end{corollary}
\state Proof. Following the proof of Theorem \ref{tConstentApproxComp}, we see that the present assumption suffices for $\nliminf \phi^\nu(x^\nu)$ $\geq$ $\phi(x)$ to hold whenever $x^\nu\to x$. By \cite[Prop. 7.7]{VaAn}, the assertion about level-sets follows. Again following the proof, we deduce that $\nOutLim (\gph S^\nu) \subset \gph S$ holds without the conditions (a)-(e) in the theorem and we invoke \cite[Thm. 5.37]{VaAn} to reach the conclusion.\eop

The assumptions of the corollary permit the approximating $\phi^\nu$ to be arbitrarily ``high'' relative to $\phi$, i.e., $\nOutLim (\epi \phi^\nu)$ could be a strict subset of $\epi \phi$. This is a main reason why we avoid calling the approximations consistent in this case. Still, the mild assumptions make the corollary widely applicable.

\section{Examples}

An array of examples illustrate the breadth of the framework, but numerous possibilities are omitted including those involving polyhedral approximations of $X$ as in \cite{BenTalNemirovski.01}. We recall that $X,X^\nu$ and $f_i, f_i^\nu$ are nonconvex unless specified otherwise.

\begin{example}{\rm (goal optimization).}\label{eGoal}
For parameters $\alpha_i\in [0, \infty)$ and $\tau_i\in \reals$, we consider the problem
\[
\nnmin_{x\in X} \nsum_{i=1}^m \alpha_i \max\big\{0,~ f_i(x) - \tau_i \big\},
\]
which aims to lower each $f_i(x)$ down to the goal of $\tau_i$, with $\alpha_i$ being the per-unit penalty for failing to do so. The problem is of the form \eqref{eqn:actualpr} with $h(z) = \nsum_{i=1}^m \alpha_i \max\{0, z_i - \tau_i \}$. A possible approximation is to set $X^\nu = X$ and $F^\nu = F$, but replace $h$ by
\[
h^\nu(z) = \nsum_{i=1}^m \alpha_i \psi^\nu(z_i - \tau_i), ~\mbox{ where } \psi^\nu(\gamma) = \frac{1}{\theta^\nu}\ln\big(1+\exp(\theta^\nu\gamma)\big), ~~~\theta^\nu\in (0,\infty).
\]
Since $h^\nu$ is convex and differentiable any number of times, it offers computational benefits over $h$. The resulting approximating problems and their optimality conditions are consistent approximations provided that $\theta^\nu\to \infty$. We refer to \cite{Chen.12,BurkeHoheiselKanzow.13} for more general smoothing schemes.
\end{example}
\state Detail. By Theorem \ref{tConstentApproxComp}, condition (a), and the refinements, we only need to show $h^\nu\eto h$. Since $\sup_{\gamma \in\reals}|\psi^\nu(\gamma)-\max\{0,\gamma\}| \leq (\ln 2)/\theta^\nu$ as seen from \cite{PolakRoysetWomersley.03}, \eqref{eqn:liminfcondition} and \eqref{eqn:limsupcondition} hold.

For an illustration of the optimality conditions, suppose that $f_1, \dots, f_m$ are smooth.
Then, $0 \in S^\nu(x,y,z)$ simplifies to $F(x) = z$, $y = \nabla h^\nu(z)$, $-\nabla F(x)^\top y \in N_X(x)$. In fact, this can be written as $-\nabla F(x)^\top\nabla h^\nu(F(x)) \in N_X(x)$, with
\[
\nabla h^\nu(z) = \begin{bmatrix}
  \alpha_1 \nabla \psi^\nu(z_1 - \tau_1)\\
  \vdots\\
  \alpha_m \nabla \psi^\nu(z_m - \tau_m)
\end{bmatrix}
~~~\mbox{ and }~~~~ \nabla \psi^\nu(\gamma) = \frac{\exp(\theta^\nu \gamma)}{1+\exp(\theta^\nu \gamma)}.
\]
Similarly, the optimality condition $0 \in S(x,y,z)$ specializes to $F(x) = z$, $y \in \partial h(z)$, $-\nabla F(x)^\top y \in N_X(x)$. Since $\partial h(z) = C_1 \times \cdots \times C_m$,  where $C_i = \{0\}$ if $z_i < \tau_i$, $C_i = [0,\alpha_i]$ if $z_i = \tau_i$, and $C_i = \{\alpha_i\}$  otherwise, the optimality condition simplifies further to $-\nabla F(x)^\top y \in N_X(x)$ with $y_i = 0$ if $f_i(x) < \tau_i$, $y_i \in  [0,\alpha_i]$  if $f_i(x) = \tau_i$, and $y_i = \alpha_i$  otherwise. These conditions are necessary for optimality in the actual and approximating problems by Proposition \ref{prop:optimality}; the qualification \eqref{eqn:qualGenOptCondComposite} holds in each case because $\dom h = \dom h^\nu = \reals^m$.\eop

\begin{example}{\rm (stochastic optimization).}\label{eStoch}
For probabilities $p = (p_1, \dots, p_m) \in P$, where $P\subset\reals^m$ is the set of nonnegative vectors with components summing to one, consider the problem
\[
\nnmin_{x\in X} \nsum_{i=1}^m p_i f_i(x),
\]
which appears in stochastic optimization and machine learning. The problem is of the form \eqref{eqn:actualpr} with $h(z) = \nsum_{i=1}^m p_i z_i$. In practice, the probabilities might not be fully known and this leads to approximating problems with $p$ replaced by $p^\nu\in P$, i.e., $h^\nu(z) = \nsum_{i=1}^m p_i^\nu z_i$ in \eqref{eqn:approxpr}.  The change may also be part of a sensitivity analysis such as when developing influence functions \cite{KohLiang.17,BasuPopeFeizi.20}. The approximations are consistent provided that $p^\nu\to p$.
\end{example}
\state Detail. Since $h^\nu(z^\nu) \to h(z)$ whenever $z^\nu\to z$, we obtain $h^\nu\eto h$; see \eqref{eqn:liminfcondition} and \eqref{eqn:limsupcondition}. Thus, the consistency follows via condition (a) of Theorem \ref{tConstentApproxComp}. The optimality condition $0 \in S(x,y,z)$ is necessary by Proposition \ref{prop:optimality} because the qualification \eqref{eqn:qualGenOptCondComposite} holds. In fact, the optimality condition simplifies to $0 \in \nsum_{i=1}^m p_i \con\partial f_i(x) + N_X(x)$ because $y = \nabla h(z) = p$.\eop

\begin{example}{\rm (distributionally robust optimization).}\label{eDistr}
Let $P\subset\reals^m$ be as in the previous example and let $A$ and $A^\nu$ be nonempty closed subsets of $P$. We consider the problem
\[
\nnmin_{x\in X} \,\nmax_{p \in A} \nsum_{i=1}^m p_i f_i(x)
\]
and its approximation obtained by replacing $A$ by $A^\nu$. The set $A$ might consist of a single, true probability vector and be approximated by a set $A^\nu$ ``centered'' on a current best estimate of the true probability vector with some ``radius'' reflecting the uncertainty about this estimate. The resulting problems fit the forms of \eqref{eqn:actualpr} and \eqref{eqn:approxpr} with $h(z) = \nmax_{p \in A} \langle p, z\rangle$  and $h^\nu(z) = \nmax_{p \in A^\nu} \langle p, z\rangle$. If $A^\nu \sto A$, then the approximations are consistent.
\end{example}
\state Detail. Let $\hat A(z) = \nargmax_{p \in A} \langle p, z\rangle$ and $\hat A^\nu(z) = \nargmax_{p \in A^\nu} \langle p, z\rangle$. The functions $h$ and $h^\nu$ are convex. They are also real-valued, and thus continuous, because $A$ and $A^\nu$ are compact. To establish $h^\nu\eto h$, we argue using \eqref{eqn:liminfcondition} and \eqref{eqn:limsupcondition}. Suppose that $z^\nu \to z$. Let $\bar p\in \hat A(z)$. Since $A^\nu \sto A$, there exists $p^\nu \in A^\nu \to \bar p$. This implies that
\[
\nliminf \big( \nmax_{p \in A^\nu} \langle p, z^\nu\rangle \big)\geq \nliminf \langle p^\nu, z^\nu\rangle = \langle \bar p, z\rangle = \nmax_{p \in A} \langle p, z\rangle.
\]
Thus, \eqref{eqn:liminfcondition} holds for $h$ and $h^\nu$. For fixed $z$, let $p^\nu \in \hat A^\nu(z)$. Then, for any $\{\bar p^\nu\in A, \nu\in\nats\}$, one has
\[
\nmax_{p \in A^\nu} \langle p, z\rangle = \langle p^\nu, z\rangle = \langle p^\nu-\bar p^\nu, z\rangle + \langle \bar p^\nu, z\rangle \leq \|p^\nu - \bar p^\nu\|_2 \|z\|_2 + \nmax_{p \in A} \langle p, z\rangle.
\]
Since there are $\bar p^\nu\in A$ such that $\|p^\nu - \bar p^\nu\|_2\to 0$ because $A^\nu \sto A$, $\nlimsup h^\nu(z) \leq h(z)$. Consequently, $h^\nu\eto h$ and we again invoke condition (a) of Theorem \ref{tConstentApproxComp} and the refinements.\eop

\begin{example}{\rm (augmented Lagrangian methods).}\label{eAug}
With $h(z) = z_1 + \sum_{i=2}^m \iota_{\{0\}}(z_i)$, the problem
\[
\nnmin_{x\in X} f_1(x) ~\mbox{ subject to } ~f_2(x) = 0, \dots, f_m(x) = 0
\]
fits the mold \eqref{eqn:actualpr}. An approximation stemming from augmented Lagrangian methods utilizes
\[
h^\nu(z) = z_1 + \nsum_{i=2}^m \big(y_i^\nu z_i + \half \theta^\nu z_i^2\big),
\]
where $y^\nu = (y_2^\nu, \dots, y_m^\nu) \in \reals^{m-1}$ and $\theta^\nu \in [0,\infty)$. If $\{y^\nu, \nu\in \nats\}$ is bounded and $\theta^\nu \to \infty$, then the resulting approximations are consistent.
\end{example}
\state Detail. The function $h$ is not real-valued and we turn to  condition (b) of Theorem \ref{tConstentApproxComp}. For $z \in$ $\dom h$ $=$ $\reals \times \{0\}^{m-1}$, one has $h^\nu(z) = z_1 = h(z)$ and the condition holds. This also confirms the limsup-condition \eqref{eqn:limsupcondition} for $h^\nu$ and $h$. For the liminf-condition \eqref{eqn:liminfcondition}, suppose that $z^\nu\to z$. Then,
\[
\nliminf h^\nu(z^\nu) \geq z_1 + \nsum_{i=2}^m \nliminf (y_i^\nu z_i^\nu) + \half \nsum_{i=2}^m \nliminf \big(\theta^\nu (z_i^\nu)^2\big).
\]
If any $z_2, \dots, z_m$ is nonzero, then the right-hand side equals $\infty$, which coincides with $h(z)$. If $z_2= \cdots = z_m = 0$, then the right-hand size is no smaller than $z_1 = h(z)$. Thus, $h^\nu  \eto h$.\eop

\begin{example}{\rm (min-functions).}\label{eMinfcn}
For smooth $g_{ik}:\reals^n\to \reals$, $k=1, \dots, s_i$, $i=1, \dots, m$, let
\[
f_i(x) = \min_{k = 1, \dots, s_i} g_{ik}(x).
\]
We adopt smooth approximations of $f_i$ by defining
\[
f_i^\nu(x) = -\frac{1}{\theta^\nu}\ln\Big(\nsum_{k=1}^{s_i}\exp\big(-\theta^\nu g_{ik}(x)\big)\Big),
\]
where $\theta^\nu \in (0,\infty)$; see \cite{Nesterov.05,Chen.12,BeckTeboulle.12, BurkeHoheiselKanzow.13} for related smoothing techniques. If $\theta^\nu\to \infty$, then \eqref{eqn:suffFconv} holds. Moreover, if in addition $X^\nu \supset X$, $X^\nu \sto X$, with these sets being convex, $h^\nu\eto h$, $h^\nu(z) \leq h^\nu(\bar z)$ when $z\leq \bar z$, and $h^\nu (z) \leq h(z)$ for all $z$, then the approximations are weakly consistent.
\end{example}
\state Detail. We leverage condition (c) of Theorem \ref{tConstentApproxComp} and observe that $F^\nu(x) \leq F(x)$ because (see \cite{PolakRoysetWomersley.03})
\begin{equation}\label{eqn:smoothingerror}
0\leq  f_i(x) - f_i^\nu(x)   \leq \frac{\ln s_i}{\theta^\nu} ~~~\forall x\in\reals^n.
\end{equation}
Moreover, we need to confirm  \eqref{eqn:suffFconv}. Let $x^\nu\to \bar x$. The previous inequalities show that $f_i^\nu(x^\nu)\to f_i(\bar x)$ because $\theta^\nu\to \infty$. It is apparent that $f_i^\nu$ is smooth and, in fact,
\[
\nabla f_i^\nu(x) = \nsum_{k=1}^{s_i} \mu^\nu_{ik}(x)\nabla g_{ik}(x), ~~\mbox{ with }~
\mu^\nu_{ik}(x)=\frac{\exp\big(\theta^\nu(f_i(x)-g_{ik}(x))\big)}{\nsum_{j=1}^{s_i}\exp\big(\theta^\nu(f_i(x)-g_{ij}(x))\big)}.
\]
Since $\mu^\nu_{ik}(x) \in (0,1)$ regardless of $x$, we have that $\{\nabla f_i^\nu(x^\nu), \nu\in \nats\}$ is bounded. For some $N\in \cN^\grill_\infty$, suppose that $\nabla f_i^\nu(x^\nu)\Nto \bar v$. We would like to show $\bar v \in \con \partial f_i(\bar x)$.
We observe that
\[
\con \partial f_i(x) = -\con \partial (-f_i)(x) = - \con \big\{ \nabla (-g_{ik})(x) ~\big|~ k\in \bbA_i(x) \big\} = \con \big\{ \nabla g_{ik}(x) ~\big|~ k\in \bbA_i(x) \big\},
\]
where $\bbA_i(x) = \nargmin_{k=1, \dots, s_i} g_{ik}(x)$.  If $k\not\in \bbA_i(\bar x)$, then the continuity of $f_i$ and $g_{ik}$ implies that
\[
\exp\big(\theta^\nu(f_i(x^\nu)-g_{ik}(x^\nu))\big)\Nto 0.
\]
The denominator in the defining expression for $\mu^\nu_{ik}(x^\nu)$ is greater than one because $f_i(x^\nu) - g_{ij}(x^\nu) = 0$ for $j\in \bbA_i(x^\nu)$. Thus, $\mu^\nu_{ik}(x^\nu)\Nto 0$.

For any $k$, $\mu^\nu_{ik}(x^\nu) \in (0,1)$. Consequently, after passing to another subsequence which we also denote by $N$, there are $\mu^\infty_{ik}\in [0,1]$, $k=1, \dots, s_i$, such that $\mu^\nu_{ik}(x^\nu) \Nto \mu_{ik}^\infty$. Since $\nsum_{k=1}^{s_i} \mu^\nu_{ik}(x^\nu) = 1$ for all $\nu$, we must also have $\nsum_{k=1}^{s_i} \mu_{ik}^\infty = 1$, with $\mu_{ik}^\infty = 0$ if $k\not\in \bbA_i(\bar x)$ as already seen. We conclude that
\[
\nabla f^\nu_i(x^\nu) = \nsum_{k=1}^{s_i} \mu^\nu_{ik}(x^\nu)\nabla g_{ik}(x^\nu) \Nto \bar v = \nsum_{k \in \bbA_i(\bar x)} \mu^\infty_{ik} \nabla g_{ik}(\bar x) \in \con \partial f_i(\bar x)
\]
and \eqref{eqn:suffFconv} holds. An instance of $h$ and $h^\nu$ satisfying condition (c) of Theorem \ref{tConstentApproxComp} follows next.\eop

\begin{example}{\rm (penalty methods).}\label{ePen}
With $h(z) = z_1 + \iota_{(-\infty,0]^{m-1}}(z_2, \dots, z_m)$, we consider the problem
\[
\nnmin_{x\in X} f_1(x) \mbox{ subject to } f_2(x) \leq 0, \dots, f_m(x) \leq 0,
\]
which is of the form \eqref{eqn:actualpr}. An approximation stemming from penalty methods utilizes
\[
h^\nu(z) = z_1 + \theta^\nu \nsum_{i=2}^m \big(\max\{0, z_i\}\big)^2,
\]
where $\theta^\nu \in [0,\infty)$. If $X^\nu = X$, $F^\nu = F$, and $\theta^\nu \to \infty$, then the resulting approximations are consistent.
\end{example}
\state Detail. We turn to condition (c) of Theorem \ref{tConstentApproxComp}. Certainly, $h^\nu(z) \leq h^\nu(\bar z)$ when $z\leq \bar z$ and the other requirements hold as well. To confirm $h^\nu\eto h$, we leverage the characterization \eqref{eqn:liminfcondition} and \eqref{eqn:limsupcondition}. Let $z^\nu\to z$. We note that $\nliminf h^\nu(z^\nu) \geq h(z)$ when $z=(z_1, \dots, z_m)$ has $z_i \leq 0$ for $i=2, \dots, m$. If any of these $z_i$ are positive, then $z_i^\nu>0$ for sufficiently large $\nu$. Thus, $h^\nu(z^\nu)\to \infty = h(z)$ and \eqref{eqn:liminfcondition} holds for $h^\nu$ and $h$. For any $z$, $\nlimsup h^\nu(z) \leq h(z)$ and \eqref{eqn:limsupcondition} holds as well. Consequently, $h^\nu  \eto h$.

Condition (c) of Theorem \ref{tConstentApproxComp} also permits outer approximations of $X$ and lower bounding approximations of $F$ (see Example \ref{eMinfcn}) while retaining weak consistency.\eop

\begin{example}{\rm (interior-point methods).}\label{eInterior}
For the actual problem in Example \ref{ePen}, we consider the approximation
\[
h^\nu(z) = \begin{cases}
z_1 - \frac{1}{\theta^\nu}\nsum_{i=2}^m \ln (- z_i) & \mbox{ if } z_2<0, \dots, z_m < 0\\
\infty & \mbox{ otherwise,}
\end{cases}
\]
where $\theta^\nu \in (0,\infty)$. This logarithmic penalty approach is the basis for (primal) interior-point methods. If $\theta^\nu \to \infty$ and for all $x\in X$, with $f_i(x) \leq 0$ for all $i=2, \dots, m$, there exists $x^\nu\in X\to x$ such that $f_i(x^\nu) < 0$ for all $i=2, \dots, m$, then the approximations are consistent.
\end{example}
\state Detail. We leverage condition (d) of Theorem \ref{tConstentApproxComp}. Trivially, $h$ is continuous relative to its domain. We establish $h^\nu \eto h$ via \eqref{eqn:liminfcondition} and \eqref{eqn:limsupcondition}. Let $z^\nu \to z$. If $z_i >0$ for some $i =2, \dots, m$, then $z_i^\nu\geq 0$ for sufficiently large $\nu$ and $\nliminf h^\nu(z^\nu) \geq h(z)$ because both sides are $\infty$. If $z_i\leq 0$ for $i=2, \dots, m$, then $h(z) = z_1$ and the liminf-condition holds again. The condition $\nlimsup h^\nu(z^\nu) \leq h(z)$ holds because for any $z_i = 0$, one can select $z_i^\nu = -\exp(-\sqrt{\theta^\nu})$. Then, $z_i^\nu \to z_i$ and $-(\ln (- z_i^\nu))/\theta^\nu = \sqrt{\theta^\nu}/\theta^\nu \to 0$.\eop

\begin{example}{\rm (expectation functions).}\label{eIntegral}
Expectation functions arise in stochastic optimization and machine learning and then the component functions $f_1, \dots, f_m$ may take the form
\[
f_i(x) = \Ex\big[g_i(\bfxi,x)\big],
\]
where $g_i:\Xi\times \reals^n\to \Reals$ is defined in terms of a probability space $(\Xi,\cB,\bbP)$. An independent and identically distributed sample $\bfxi_1, \dots, \bfxi_\nu$ according to $\bbP$ defines a sample average approximation
\[
f_i^\nu(x) = \frac{1}{\nu}\nsum_{j=1}^\nu g_i(\bfxi_j, x).
\]
If both $g_i$ and $-g_i$ are random lower semicontinuous and locally inf-integrable (see \cite[Sec. 8.G]{primer} for definitions), then, with probability one, $f_i^\nu(x^\nu) \to f_i(x)$ whenever $x^\nu\to x$ by \cite[Thm. 8.56]{primer}.
\end{example}
\state Detail. If there is a set $B\in \cB$ such that $\bbP(B) = 1$ and $g_i(\xi,\cdot)$ is smooth for all $\xi\in B$, then one can attempt to  repeat the above arguments with $f_i$ replaced by the partial derivatives of $f_i$ and conclude that $\nabla f^\nu_i(x^\nu) \to \nabla f_i(x)$ whenever $x^\nu\to x$. This would allow us to satisfy the requirement \eqref{eqn:suffFconv}. For further details and refinements about approximations of expectation functions and their subgradients, we refer to \cite{DavisDrusvyatskiy.22}.\eop

\begin{example}{\rm (oracle functions).}\label{eOuterapprox}
Suppose that $h$ is real-valued, but not available in an explicit form. If for each $z$ we can compute $h(z)$ and a subgradient $v\in \partial h(z)$, then the approximation
\[
h^\nu(z) = \max_{k=1, \dots, \nu} h(z^k) + \langle v^k, z - z^k\rangle,~~\mbox{ with } ~z^k\in\reals^m, ~~v^k \in \partial h(z^k),
\]
remains available. Now, $h^\nu \eto h$ if $\{z^k, k\in\nats\}$ is a countable dense subset of $\reals^m$.
\end{example}
\state Detail. By \cite[Thm. 7.17]{VaAn}, $h^\nu \eto h$ whenever the functions converge pointwise on a countable dense subset of $\reals^m$. Let $\bar z\in \{z^k, k\in\nats\}$. Obviously, $h^\nu(\bar z) \leq h(\bar z)$ because $h$ is the pointwise supremum of its affine supports; see \cite[Thm. 8.13]{VaAn}. Moreover, there is $\bar \nu$ such that $\bar z\in \{z^1, \dots, z^{\bar\nu}\}$. Thus, $h^\nu(\bar z) \geq h(\bar z)$ for all $\nu\geq \bar \nu$. By \cite[Thm. 7.17]{VaAn}, this also means that $h^\nu(z)\to h(z)$ for all $z\in\reals^m$.\eop

\begin{example}{\rm (homotopy method).}\label{eHomotopy}
For proper, lsc, and convex $\hat h:\reals^{m-1}\to \Reals$ and lLc $\hat F:\reals^n \to \reals^{m-1}$, consider the problem
\[
\nnmin_{x\in X} \hat h\big( \hat F(x)\big),
\]
which is of the form \eqref{eqn:actualpr} with $h(z) = \hat h(z_1, \dots, z_{m-1})$ for $z = (z_1, \dots, z_m)$ and $F(x) = (\hat F(x), f_m(x))$ for some lLc $f_m:\reals^n\to \reals$. A homotopy method for the problem solves the approximations
\[
\nnmin_{x\in X} \,(1-\lambda^\nu) \hat h\big(\hat F(x)\big) + \lambda^\nu f_m(x)
\]
as $\lambda^\nu \to 0$. If $f_m$ is chosen wisely, then the approximating problems might be simpler to solve than the actual one while benefitting from warm starts. The approximations are of the form \eqref{eqn:approxpr} with
\[
h^\nu(z) = (1-\lambda^\nu) \hat h(z_1, \dots, z_{m-1}) + \lambda^\nu z_m.
\]
If $\lambda^\nu\to 0$, then the approximations are consistent.
\end{example}
\state Detail. To establish $h^\nu\eto h$ via \eqref{eqn:liminfcondition} and \eqref{eqn:limsupcondition}, we note that for $z^\nu\to z$, one has
\[
\nliminf (1-\lambda^\nu) \hat h(z_1^\nu, \dots, z_{m-1}^\nu) + \lambda^\nu z_m^\nu \geq  \nliminf \hat h(z_1^\nu, \dots, z_{m-1}^\nu) + \nliminf \lambda^\nu z_m^\nu \geq \hat h(z_1, \dots, z_{m-1}) = h(z).
\]
Moreover, for $z\in \dom h$, $(1-\lambda^\nu) \hat h(z_1, \dots, z_{m-1}) + \lambda^\nu z_m  \to  \hat h(z_1, \dots, z_{m-1})$. We therefore have both epi-convergence and pointwise convergence (on $\dom h$) and condition (b) of Theorem \ref{tConstentApproxComp}, with refinements, establishes consistency.\eop

\begin{example}{\rm (monitoring functions).}\label{eMonitoring}
In extended nonlinear programming \cite{Rockafellar.99}, one utilizes
\[
h(z) = \nsup_{y\in Y}\lset    \langle z,y \rangle  - \half\langle y,By \rangle\rset ~~~\mbox{ and }~~~ h^\nu(z) = \nsup_{y\in Y^\nu}\lset    \langle z,y \rangle  - \half\langle y,B^\nu y \rangle\rset,
\]
where $Y,Y^\nu\subset\reals^m$ are nonempty polyhedral sets and  $B,B^\nu$ are symmetric positive semidefinite $m\times m$-matrices.
If $B^\nu \to B$ and $Y^\nu \sto Y$, then $h^\nu\eto h$. If in addition $B$ is positive definite or $Y$ is bounded, then  $h$ is real-valued and $h^\nu(z)\to h(z)$ for $z\in \reals^m$. Thus, main steps toward (weak) consistency via Theorem \ref{tConstentApproxComp} are immediately accomplished.
\end{example}
\state Detail. The epi-convergence $h^\nu \eto h$ is established via the corresponding conjugate functions $h^{\nu*}$ and $h^*$. By Wijsman's theorem \cite[Thm. 11.34]{VaAn}, it suffices to show that $h^{\nu*}\eto h^*$. Since $h^*(y) = \half \langle y, By\rangle + \iota_Y(y)$, $h^{\nu*}(y) = \half \langle y, B^\nu y\rangle + \iota_{Y^\nu}(y)$, and $\iota_{Y^\nu}\eto \iota_Y$, we conclude that $h^{\nu*}\eto h^*$.

It follows by \cite[Thm. 7.17]{VaAn} that the pointwise convergence $h^\nu(z) \to h(z)$ holds when $h$ is real-valued. This is the case when $B$ is positive definite or when $Y$ is bounded.\eop

\begin{example}{\rm (difference-of-convex functions).}\label{eDC}
For a proper, lsc, and convex $f:\reals^n\to \Reals$ and convex $g:\reals^n\to \reals$, consider the problem
\[
\nnmin_{x\in X} f(x) - g(x),
\]
which involves an objective function of the difference-of-convex kind. As the only approximation, suppose that $g$ is replaced by smooth functions $g^\nu:\reals^n\to \reals$ that also satisfy the property:
\begin{equation*}\label{eqn:assDC}
x^\nu \in X \to x~~\Longrightarrow~~ g^\nu(x^\nu) \to g(x) \mbox{ and } \big\{\nabla g^\nu(x^\nu),  \nu\in\nats\big\} \mbox{ is bounded with all cluster points in } \partial g(x),
\end{equation*}
which holds, for instance, when $g^\nu$ is convex and $g^\nu\eto g$; see the arguments leading to the third refinement of Theorem \ref{tConstentApproxComp}. This produces weakly consistent approximations.
\end{example}
\state Detail. For $z = (z_1, z_2)$, with $z_1\in \reals^n$ and $z_2 \in \reals$,  we set $h(z) = h_1(z_1) + h_2(z_2)$, where $h_1(z_1) = f(z_1)$ and $h_2(z_2) = -z_{2}$. Moreover, let $F(x) = (x, g(x))$ and $F^\nu(x) = (x, g^\nu(x))$ so that $m = n+1$, $m_1 = n$, $m_2 = 1$, and $r = 2$ in the notation of condition (e) of Theorem \ref{tConstentApproxComp}. Since $g$ is real-valued and convex, it is lLc and then $F$ is also lLc. With $h^\nu = h$, we trivially obtain $h^\nu \eto h$. The requirement \eqref{eqn:suffFconv} translates into the assumed property. We then invoke condition (e) of Theorem \ref{tConstentApproxComp} to conclude weak consistency. For instance $g(x) = \max_{k=1, \dots, s} g_k(x)$, with $g_k:\reals^n\to \reals$ being smooth, could be approximated parallel to Example \ref{eMinfcn}. In this case, $0 \in S(x,y,z)$ reduces to $0 \in \partial f(x) - \partial g(x) + N_X(x)$.\eop

\section{Enhanced Proximal Composite Algorithm}

As an example of an implementable version of the consistent approximation algorithm, we consider the setting where $X^\nu$ is convex, $h^\nu$ is real-valued, and $F^\nu$ is twice continuously differentiable. Then, the approximating problems \eqref{eqn:approxpr} are solvable by proximal composite methods, which can be traced back to \cite{Fletcher.82,Powell.83,Burke.85}; see also \cite{Powell.84,Yuan.85,BurkeFerris.95} for trust-region versions and \cite{DrusvyatskiyPaquette.19} for details about rates of convergence under the assumption that $h^\nu$ and $\nabla F^\nu$ are Lipschitz continuous. Let $\prj_{X^\nu}(x)$ be the projection of $x$ on $X^\nu$.

\bigskip

\state Enhanced Proximal Composite Algorithm (EPCA).

\begin{description}

  \item[Data.] ~\,~$x^0\in \reals^n$, $\tau\in (1,\infty), \sigma\in (0,1)$, $\bar \lambda \in (0,\infty)$, $\lambda^0 \in (0, \bar \lambda]$, $\{\delta^\nu> 0, \nu\in\nats\}\to 0$.

  \item[Step 0.]  Set $\nu = 1$.

  \item[Step 1.]  Set $k = 0$ and $\bar x^0 = \prj_{X^\nu}(x^{\nu-1})$.

  \item[Step 2.]  Compute $x^\star \in \nargmin_{x\in X^\nu} h^\nu\big(F^\nu(\bar x^k) + \nabla F^\nu(\bar x^k)(x-\bar x^k)\big) + \tfrac{1}{2\lambda^k}\|x-\bar x^k\|_2^2$.

~\,~~~If $x^\star = \bar x^k$, then go to Step 4.

  \item[Step 3.]  If $h^\nu\big(F^\nu(\bar x^k)  \big) - h^\nu\big(F^\nu(x^\star)  \big) \geq \sigma \Big( h^\nu\big(F^\nu(\bar x^k)\big) - h^\nu\big(F^\nu(\bar x^k) + \nabla F^\nu(\bar x^k)(x^\star - \bar x^k) \big) \Big)$,

~~~~~~~~~then set $\lambda^{k+1} = \min\{\tau \lambda^k, \bar\lambda\}$ and go to Step 5.

~~~~\,Else, replace $\lambda^k$ by $\lambda^k/\tau$ and go to Step 2.

  \item[Step 4.] Set $x^\nu = x^\star$, $z^\nu = F^\nu(x^\nu)$, and $y^\nu$ such that $y^\nu \in \partial h^\nu(z^\nu)$ and $-\nabla F^\nu(x^\nu)^\top y^\nu \in N_{X^\nu}(x^\nu)$.

  ~~~~\,Replace $\nu$ by $\nu +1$ and go to  Step 1.

  \item[Step 5.] Set $\bar x^{k+1} = x^\star$, $\bar z^{k+1} = F^\nu(\bar x^k) + \nabla F^\nu(\bar x^k) (\bar x^{k+1}  - \bar x^k)$, and $\bar y^{k+1}$ such that
\[
 \bar y^{k+1} \in \partial h^\nu(\bar z^{k+1}) ~~\mbox{ and }~~ -\nabla F^\nu(\bar x^k)^\top \bar y^{k+1} - \tfrac{1}{\lambda^k}(\bar x^{k+1} - \bar x^k)
  \in N_{X^\nu}(\bar x^{k+1}).
\]
~~~~\,Set $\bar u^{k+1}  = F^\nu(\bar x^{k+1}) - \bar z^{k+1}$ and $\bar w^{k+1}  = \big(\nabla F^\nu(\bar x^{k+1}) - \nabla F^\nu(\bar x^k)\big)^\top \bar y^{k+1}  - \tfrac{1}{\lambda^k} (\bar x^{k+1} - \bar x^k)$.

~\,\,\,\,\,~If $\max\big\{\|\bar u^{k+1}\|_2, \|\bar w^{k+1}\|_2\big\}\leq\delta^\nu$, then set $x^\nu = \bar x^{k+1}$, $y^\nu = \bar y^{k+1}$, $z^\nu = \bar z^{k+1}$, replace $\nu$ by $\nu +1$,

~~~~~~~~~and go to  Step 1.

~\,~~~Else, replace $k$ by $k +1$, and go to  Step 2.

\end{description}

\medskip

\begin{theorem}{\rm (enhanced proximal composite algorithm).}\label{tEPCA}
Suppose that $X^\nu \sto X$ with these sets being convex, $h^\nu \eto h$ with $h^\nu$ real-valued, and \eqref{eqn:suffFconv} holds with $f_i^\nu$ twice continuously differentiable for $i=1, \dots, m$. If $\{(x^\nu,y^\nu,z^\nu), \nu\in\nats\}$ is generated by EPCA and $(x^\nu,y^\nu,z^\nu)\Nto(\hat x, \hat y, \hat z)$ for some $N\in \cN_\infty^\grill$, then
\[
\hat x\in X, ~~~~~~h\big(F(\hat x)\big) \leq \nliminf_{\nu\in N} h^\nu\big(F^\nu(x^\nu)\big), ~~~~~~ 0 \in S(\hat x,\hat y,\hat z).
\]
\end{theorem}
\state Proof. The two first claims follow by Corollary \ref{cApprox}. The construction in Steps 4 and 5 ensures that $(x^\nu,y^\nu,z^\nu) \in (S^\nu)^{-1}(\ball(0,\delta^\nu))$ for all $\nu$ when the norm on $\reals^{2m+n}$ is $(u,v,w)\mapsto \max\{\|u\|_2, \|v\|_2, \|w\|_2\}$; see the discussion below. Thus, Corollary \ref{cApprox} also confirms that $0 \in S(\hat x,\hat y,\hat z)$.\eop

EPCA solves the sequence of approximating problems using a composite proximal method of the kind proposed in \cite{LewisWright.16} for each problem. The outer loop indexed by $\nu$ represents refinement of the approximating problems. The inner loop indexed by $k$ corresponds to the iterative solution of each approximating problem by means of solving the convex subproblems in Step 2 obtained by linearization of $F^\nu$. (For inexact solution of such subproblems, we refer to \cite{DrusvyatskiyPaquette.19}.) The proximal parameter $\lambda^k$ is controlled adaptively using the test in Step 3.

A concern is whether EPCA terminates after a finite $k$ in Step 4 or 5 and thus always produces the next $(x^\nu, y^\nu, z^\nu)$. This turns out to be the case as long as
\begin{equation}\label{eqn:levelassumption}
\big\{ \iota_{X^\nu} + h^\nu \circ F^\nu \leq  h^\nu\big(F^\nu(\bar x^{0})\big)\big\} ~\mbox{ is bounded,}
\end{equation}
where $\bar x^0 = \prj_{X^\nu}(x^{\nu-1})$. We see this as follows:

An optimality condition for the subproblem in Step 2 states that $x^\star$ satisfies
\begin{equation}\label{eqn:proxalgderiv1}
-\nabla F^\nu(\bar x^k)^\top y - \tfrac{1}{\lambda^k} (x^\star - \bar x^k) \in N_{X^\nu}(x^\star) ~\mbox{ for some } ~ y \in \partial h^\nu\big( F^\nu(\bar x^k) + \nabla F^\nu(\bar x^k) (x^\star - \bar x^k) \big).
\end{equation}
Thus, when $x^\star = \bar x^k$, one has $y \in \partial h^\nu(F^\nu(\bar x^k))$ and $-\nabla F^\nu(\bar x^k)^\top y \in N_{X^\nu}(x^\star)$. This means that there exists $y^\nu$ in Step 4. The best way to compute such $y^\nu$ depends on the nature of $h^\nu$ and $X^\nu$, but it can be achieved by convex optimization since it at most involves finding points in two convex sets. Then, $(x^\nu,y^\nu,z^\nu) \in (S^\nu)^{-1}(0)$. In Step 5, there exists likewise $\bar y^{k+1}$ by \eqref{eqn:proxalgderiv1}, again computable by convex optimization. The construction in Step 5 ensures that $(\bar u^{k+1}, 0, \bar w^{k+1}) \in S^\nu(\bar x^{k+1}, \bar y^{k+1}, \bar z^{k+1})$.

For the sake of contradiction, suppose that the algorithm iterates indefinitely without generating the next $(x^\nu,y^\nu,z^\nu)$. Let $\{\bar x^k, k\in\nats\}$ be the resulting sequence, which has a cluster point in view of \eqref{eqn:levelassumption}, say $\bar x$. We deduce from \cite{LewisWright.16}, especially Theorem 5.4 and is proof, that $\bar x^{k+1} - \bar x^k\to 0$ as $k\to \infty$ and $(\bar x^{k+1}-\bar x^k)/\lambda^k\to 0$ as $k\to \infty$ along the subsequence corresponding to $\bar x$. (We note that EPCA eventually exits Steps 2-3 and moves to Step 4 or 5, and that $\{1/\lambda^k, k\in \nats\}$ is bounded along the subsequence corresponding to the cluster point $\bar x$.) Then, for $k$ along the subsequence corresponding to $\bar x$, $\bar z^{k+1}\to F^\nu(\bar x)$, which means that the associated sets of subgradients $\partial h^\nu(\bar z^{k+1})$ are contained in a compact set and then the same holds for $\bar y^{k+1}$. This fact as well as the recognition that $F^\nu(\bar x^{k+1}) \to F^\nu(\bar x)$  and $\nabla F^\nu(\bar x^{k+1})-\nabla F^\nu(\bar x^{k}) \to 0$ imply that $\bar u^{k+1}\to 0$ and $\bar w^{k+1}\to 0$ as $k\to\infty$ along the subsequence. Consequently, there's $k^\star$ such that $\max\{\|\bar u^{k^\star+1}\|_2,  \|\bar w^{k^\star+1}\|_2\} \leq \delta^\nu$. We have shown that EPCA augments $\nu$ after some finite $k$ in either Step 4 or  5.

\section{Inverse Problems in Machine Learning}

An inverse problem in machine learning is that of determining an input to a collection of neural networks such that their outputs best match a given quantity \cite{HandVoroninski.18,HuangHandHeckelVoroninski.21}. Specifically, we are given $s$ neural networks represented by the mappings $F_i:\reals^{n_0}\to \reals^{n_q}$, $i=1, \dots, s$. Each network associates an input vector $x \in \reals^{n_0}$ with an output vector $F_i(x)\in \reals^{n_q}$. The goal is to determine an input vector $x$ such that the output vectors $F_1(x), \dots, F_s(x)$ are optimized in the sense of a convex function $\hat h:\reals^{s n_q}\to \reals$, i.e.,
\begin{equation}\label{eqn:actualMLproblem}
\nnmin_{x\in \hat X} ~\hat h\big(F_1(x), \dots, F_s(x)\big),
\end{equation}
where $\hat X\subset \reals^{n_0}$ is a nonempty closed set that could impose restrictions on the choice of input vector. For some target output vector $t\in  \reals^{n_q}$ and nonnegative weights $p_1, \dots, p_s$, we may simply have
\[
\hat h\big(F_1(x), \dots, F_s(x)\big) =  \nsum_{i=1}^s p_i \big\| t - F_i(x)\big\|_2^2.
\]

Suppose that the neural networks are of the feed-forward kind producing their output from passing through $q$ affine mappings composed with an activation function. For each layer $k=1, \dots, q$, there is a given $n_k \times n_{k-1}$-matrix $A_{i,k}$, a given vector $b_{i,k}\in \reals^{n_k}$, and $F_{i,k}:\reals^{n_{k-1}} \to \reals^{n_k}$ with
\[
F_{i,k}(x) = G_{i,k}(A_{i,k} x + b_{i,k})~~ \mbox{ for }~ x\in \reals^{n_{k-1}},
\]
where $G_{i,k}:\reals^{n_k}\to \reals^{n_k}$ has $G_{i,k}(y) = ( g_{i,1,k}(y_1), \dots, g_{i,n_k,k}(y_{n_k}) )$ for $y = (y_1, \dots, y_{n_k}) \in \reals^{n_k}$ and $g_{i,j,k}:\reals\to \reals$ are given epi-regular lLc functions. For example, the ReLU function $g_{i,j,k}(\gamma) = \max\{0, \gamma\}$ is a common such activation function; it is epi-regular and lLc. Then, $F_i = F_{i,q} \circ \dots \circ F_{i,1}$.

It is convenient to express \eqref{eqn:actualMLproblem} using additional variables as follows: We think of $x_{1,i,k} \in\reals^{n_k}$ as the output of the $k$th layer for neural network $i$. Let $r = \sum_{k=1}^q n_k$,  $X = \hat X \times \reals^{sr}$, and $x = (x_0, x_1)$, with  $x_0 \in \reals^{n_0}$ and $x_1 = (x_{1,1,1}, \dots, x_{1,s,q})\in \reals^{sr}$. Set
\[
F(x) = \big( x_{1,1,q}, \dots, x_{1,s,q},  H_1(x), \dots, H_s(x) \big) \in \reals^{s n_q + sr}~~~\mbox{ and }~~~h(z) = \hat h(z_0) + \iota_{\{0\}^{sr}}( z_1, \dots, z_s  ),
\]
where, for $z = (z_0, z_1, \dots, z_s)$, with $z_0 \in \reals^{s n_q}$ and $z_i \in \reals^{r}$, $i=1, \dots, s$, and
\[
H_i(x) = \begin{bmatrix}
  G_{i,1}\big( A_{i,1} \,x_0 + b_{i,1} \big) - x_{1,i,1}\\
  G_{i,2}\big( A_{i,2} \,x_{1,i,1} + b_{i,2} \big) - x_{1,i,2}\\
  \vdots\\
  G_{i,q}\big( A_{i,q} \,x_{1,i,q-1} + b_{i,q} \big) - x_{1,i,q}
  \end{bmatrix}.
\]
Thus, the inverse problem \eqref{eqn:actualMLproblem} is equivalently expressed in the form \eqref{eqn:actualpr} with $n = n_0 + sr$ and $m = s n_q + sr$; $h$ is proper, lsc, and convex and $F$ is lLc.
For theoretical and computational reasons, several approximations may arise. Suppose that $A_{i,k}^\nu$, $b^\nu_{i,k}$, convex $\hat h^\nu: \reals^{sn_q}\to \reals$, and epi-regular lLc $g_{i,j,k}^\nu:\reals\to \reals$  approximate $A_{i,k}$, $b_{i,k}$, $\hat h$, and $g_{i,j,k}$, respectively. The approximating quantities define $h^\nu$ and $F^\nu$ via $H_i^\nu$ and $G_{i,k}^\nu$ in the manner laid out for $h$, $F$, $H_i$, and $G_{i,k}$. For simplicity, $\hat X$ is not approximated so that $X^\nu=X$. It turns out that weak consistency follows naturally.

\begin{proposition}{\rm (weak consistency in inverse machine learning).}
In the notation of this section, suppose that for each $(i,j,k)$, the following property holds:
\begin{equation}\label{eqn:InverseProp}
  \begin{rcases}
    \gamma^\nu \in \reals \to \gamma\\
    \alpha^\nu \in \partial g_{i,j,k}^\nu(\gamma^\nu)
  \end{rcases}
  ~\Longrightarrow ~
  \begin{cases}
    g_{i,j,k}^\nu(\gamma^\nu) \to g_{i,j,k}(\gamma)\\
    \{\alpha^\nu, \nu\in\nats\} \mbox{ is bounded with all its cluster points in } \partial g_{i,j,k}(\gamma).
  \end{cases}
\end{equation}
If $A_{i,k}^\nu\to A_{i,k}$, $b^\nu_{i,k} \to b_{i,k}$, and $\hat h^\nu \eto \hat h$, then the resulting pairs $\{(\phi^\nu,S^\nu), \nu\in\nats\}$ are weakly consistent approximations of $(\phi,S)$.
\end{proposition}
\state Proof. We establish $\phi^\nu\eto \phi$ using \eqref{eqn:liminfcondition} and \eqref{eqn:limsupcondition}. The liminf-condition holds because $h^\nu \eto h$ and $H_i^\nu(x^\nu)\to H_i(x)$ whenever $x^\nu\to x$. For the limsup-condition, it suffices to consider $\bar x = (\bar x_0, \bar x_1)$ such that $\bar x_0 \in \hat X$ and $H_i(\bar x) = 0$ for $i=1, \dots, s$. Construct $x^\nu = (x_0^\nu, x_1^\nu)$ with $x_0^\nu = \bar x_0$ and $x_1^\nu = (x_{1,1,1}^\nu, \dots, x_{1,s,q}^\nu)$, where for each $i=1, \dots, s$,
\[
x_{1,i,1}^\nu = G_{i,1}^\nu\big( A_{i,1}^\nu\, \bar x_0 + b_{i,1}^\nu \big), ~~ x_{1,i,2}^\nu = G_{i,2}^\nu\big( A_{i,2}^\nu\, x_{1,i,1}^\nu + b_{i,2}^\nu \big), ~\dots, ~   x_{1,i,q}^\nu = G_{i,q}^\nu\big( A_{i,q}^\nu \,x_{1,i,q-1}^\nu + b_{i,q}^\nu \big).
\]
Thus, we have $H_i^\nu(x^\nu) = 0$ and $x^\nu \to \bar x$. Since $\hat h$ is real-valued and convex, $\hat h^\nu(x_{1,1,q}^\nu, \dots, x_{1,s,q}^\nu)$ $\to$ $\hat h(\bar x_{1,1,q}, \dots, \bar x_{1,s,q})$ by \cite[Thm. 7.17]{VaAn}. This implies that $\phi^\nu(x^\nu) \to \phi(\bar x)$.

Next, we turn to the optimality conditions. Suppose that $(\bar x, \bar y, \bar z, \bar u, \bar v, \bar w)$$\in$$\nOutLim (\gph S^\nu)$. Then, there are $N\in \cN_\infty^\grill$ and $(x^\nu, y^\nu, z^\nu, u^\nu, v^\nu, w^\nu)$$\Nto$$(\bar x, \bar y, \bar z, \bar u, \bar v, \bar w)$ with $(u^\nu,v^\nu,w^\nu)\in S^\nu(x^\nu,y^\nu,z^\nu)$.

Let $A_{i,j,k}$ be the $j$th row of $A_{i,k}$ and $b_{i,j,k}$ be the $j$th component of $b_{i,k}$, with parallel definitions of $A^\nu_{i,j,k}$ and $b^\nu_{i,j,k}$. The $j$th component of $x_{1,i,k} \in \reals^{n_k}$ is denoted by $x_{1,i,j,k}$. Moreover, the affine function $a_{i,j,k}:\reals^{n_0 + sr}\to \reals$ has
\[
a_{i,j,k}(x) = \langle A_{i,j,k},\, x_{1,i,k-1}\rangle + b_{i,j,k}
\]
with the convention that $x_{1,i,0} = x_0$. The function $f_{i,j,k}:\reals^{n_0 + sr}\to \reals$ has
\[
f_{i,j,k}(x) = g_{i,j,k}\big( a_{i,j,k}(x) \big) - x_{1,i,j,k}.
\]
The approximations $a^\nu_{i,j,k}$ and $f^\nu_{i,j,k}$ are defined similarly.

Since $u^\nu = F^\nu(x^\nu) - z^\nu$, we conclude that $\bar u = F(\bar x)-\bar z$; recall that $g_{i,j,k}^\nu(\gamma^\nu) \to g_{i,j,k}(\gamma)$ when $\gamma^\nu\to \gamma$.
The fact that $v^\nu + y^\nu \in \partial h^\nu(z^\nu)$ implies that $v^\nu_0 + y^\nu_0 \in \partial \hat h^\nu(z_0^\nu)$  and $z_i^\nu = 0$, $v_i^\nu+y_i^\nu \in\reals^r$, $i=1, \dots, s$, where $v^\nu = (v_0^\nu, v_1^\nu, \dots, v_s^\nu)$, with $v_0^\nu \in \reals^{sn_q}$ and $v_i^\nu \in \reals^r$ for $i=1, \dots, s$; the vectors $\bar v$, $\bar y$, $\bar z$, $y^\nu$, and $z^\nu$ are partitioned similarly. By Attouch's theorem \cite[Thm. 12.35]{VaAn}, we conclude that $\bar v_0 + \bar y_0 \in \partial \hat h(\bar z_0)$. Thus, $\bar v + \bar y \in \partial h(\bar z)$.

It remains to confirm the last portion of $(\bar u, \bar v, \bar w) \in S(\bar x, \bar y, \bar z)$ involving subgradients of the component functions of $F$. Since $g_{i,j,k}$ is lLc and epi-regular, the component function $f_{i,j,k}$ has
\[
\partial f_{i,j,k}(x) = c_{i,j,k}^\top \partial g_{i,j,k}\big( a_{i,j,k}(x) \big) - e_{i,j,k},
\]
where $c_{i,j,k} = \nabla a_{i,j,k}(x)$ and $e_{i,j,k}\in \reals^{n_0 + sr}$ is a vector with a single 1 placed appropriately and with zero elsewhere. Likewise,
\[
\partial f^\nu_{i,j,k}(x) = (c^\nu_{i,j,k})^\top \partial g^\nu_{i,j,k}\big( a^\nu_{i,j,k}(x)  \big) - e_{i,j,k},
\]
where $c^\nu_{i,j,k} = \nabla a^\nu_{i,j,k}(x)$. The assumption \eqref{eqn:InverseProp} implies that any sequence $\{d_{i,j,k}^\nu \in \partial f^\nu_{i,j,k}(x^\nu), \nu\in N\}$ is bounded and thus has a cluster point, with any such point in $\partial f_{i,j,k}(\bar x)$. Since $N_{\hat X}$ is outer semicontinuous, we deduce that the last requirement holds and, thus, $(\bar u, \bar v, \bar w) \in S(\bar x, \bar y, \bar z)$.\eop

A typical activation function is the ReLU  with $g_{i,j,k}(\gamma) = \max\{0, \gamma\}$. A smooth approximation could be $g^\nu_{i,j,k}(\gamma) = \frac{1}{\theta^\nu} \ln ( 1 + \exp(\theta^\nu\gamma))$. The assumption \eqref{eqn:InverseProp} holds in this case; see Example \ref{eGoal}.

\section{Rates and Error Estimates}

Consistency furnishes guarantees about the limiting behavior of approximations, but it also can be beneficial to quantify the rate of convergence. In this section, we refine results from \cite{Royset.20b} and estimate the discrepancy between near-solutions of the optimality condition $0 \in S^\nu(x,y,z)$ and those of $0\in S(x,y,z)$. Chapter 8 of \cite{CuiPang.21} addresses similar issues using different techniques, especially for problems with affine structure. For error estimates of minimizers, minima, and level-sets, we refer to \cite{Royset.18,Royset.20b}.

The {\em point-to-set distance} between $\bar x\in \reals^n$ and $C\subset \reals^n$ under norm $\|\cdot\|$ is denoted by
\[
\dist(\bar x,C) = \ninf_{x\in C} \|x-\bar x\| ~\mbox{ when }~ C\neq \emptyset~~ \mbox{ and }~ \dist(\bar x,\emptyset) = \infty.
\]
The {\em excess} of $C\subset\reals^n$ over $D\subset\reals^n$ is defined as
\[
\exs(C; D) = \begin{cases}
  \nsup_{x\in C} \dist(x, D) & \mbox{ if } C\neq \emptyset, D\neq \emptyset\\
  \infty & \mbox{ if } C\neq \emptyset, D = \emptyset\\
  0 & \mbox{ otherwise.}
\end{cases}
\]
We concentrate on a truncated version given by
\[
\exs_\rho(C; D) = \exs\big(C\cap\ball(0,\rho); D\big)~ ~\mbox{ for } \rho \in [0,\infty).
\]

All these concepts rely on the choice of norm. The Euclidean norm remains the default, but in the context of $S,S^\nu:\reals^{n+2m}\tto \reals^{2m+n}$ we adopt the norms given by
\begin{equation}\label{eqn:normCompThm}
\|(x,y,z)\|_{\rm in} = \max\big\{\|x\|_2,\|y\|_2,\|z\|_2\big\} ~~\mbox{ and } ~~\|(u,v,w)\|_{\rm out} =
\max\big\{\|u\|_2,\|v\|_2, \|w\|_2\big\}
\end{equation}
for the argument and value spaces, respectively, where $x\in\reals^n$, $y\in\reals^m$, $z\in\reals^m$, $u\in\reals^m$, $v\in\reals^m$, and $w\in\reals^n$. The graphs of $S,S^\nu$ are subsets of $\reals^{n+2m}\times \reals^{2m+n}$, for which we adopt the norm given by
\begin{equation}\label{eqn:normCompThm2}
\max\big\{\|(x,y,z)\|_{\rm in}, \|(u,v,w)\|_{\rm out}\big\}.
\end{equation}
Thus, $S^{-1}(\ball(0,\epsilon))$ is the set of near-solutions of $0\in S(x,y,z)$ with the tolerance now being specified by $\|\cdot\|_{\rm out}$, i.e., $(\bar x, \bar y, \bar z) \in S^{-1}(\ball(0,\epsilon))$ if and only if $(\bar u, \bar v, \bar w) \in S(\bar x, \bar y, \bar z)$ and $\|(\bar u, \bar v, \bar w)\|_{\rm out} \leq \epsilon$.

\begin{proposition}\label{thm:approxgeneralizedequations}{\rm (solution error in optimality conditions).} Suppose that $0 \leq\delta^\nu\leq \rho< \infty$ and $\epsilon \geq \delta^\nu + \exs_\rho(\gph S^\nu;~ \gph S)$. Then, under the norms \eqref{eqn:normCompThm} and \eqref{eqn:normCompThm2}, one has
\[
\exs_\rho\Big(\,(S^\nu)^{-1}\big(\ball(0, \delta^\nu)\big); \,~S^{-1}\big(\ball(0, \epsilon)\big)    \,\Big) \leq \exs_\rho\big(\gph S^\nu;~ \gph S\big).
\]
\end{proposition}
\state Proof. A slight modification of the proof of Theorem 5.1 in \cite{Royset.20b} yields this fact.\eop

The proposition ensures that if $(\hat x, \hat y, \hat z) \in \ball(0,\rho)$ nearly satisfies the optimality condition $0\in S^\nu(x,y,z)$ with tolerance $\delta^\nu \in [0,\infty)$, i.e., $\dist(0, S^\nu(\hat x, \hat y, \hat z)) \leq \delta^\nu$, then $(\hat x, \hat y, \hat z)$ is no further away than $\exs_\rho( \gph S^\nu;~\gph S)$ from a point that nearly satisfies the optimality condition $0\in S(x,y,z)$ with tolerance $\epsilon$. Thus, the rate of convergence of near-solutions of the approximating problems is governed by the rate of decay of $\exs_\rho( \gph S^\nu;~\gph S)$ as $\nu\to \infty$.

In the context of the consistent approximation algorithm, Proposition \ref{thm:approxgeneralizedequations} asserts that the solution tolerance $\epsilon$ (for the actual problem) is bounded by the sum of $\delta^\nu$ (the tolerance adopted when solving the approximating problem) and the error in the approximation as expressed by $\exs_\rho (\gph S^\nu; \gph S)$. Thus, it would be ``optimal'' to make $\delta^\nu$ in the consistent approximation algorithm vanish at the same rate as $\exs_\rho (\gph S^\nu; \gph S)$ tends to zero. Concrete illustrations follow in the examples below.

\begin{theorem}{\rm (estimate of excess).}\label{tConstentApproxCompQuant}
Suppose that $\rho\in [0,\infty)$ and $X= X^\nu$. Let
\begin{align*}
\eta_0^\nu & = \sup\Big\{ \big\| F^\nu(x) - F(x)\big\|_2~\Big|~x \in X \cap \ball(0,\rho)\Big\}\\
\eta^\nu & = \sup\Big\{\exs\big(\partial f_i^\nu(x); \con \partial f_i(x)\big)~\Big|~x \in X \cap \ball(0,\rho), ~i=1, \dots, m\Big\}.
\end{align*}
If the norm on $\reals^{n+2m}\times\reals^{2m+n}$ is \eqref{eqn:normCompThm2} and the norm on $\reals^m\times\reals^m$ is $\max\{\|z\|_2, \|v\|_2\}$, then one has
\[
\exs_\rho(\gph S^\nu;~ \gph S) \leq \max\big\{\sqrt{m}\rho\eta^\nu,~ \eta^\nu_0 + \exs_{2\rho}(\gph \partial h^\nu;~ \gph \partial h)\big\}.
\]
\end{theorem}
\state Proof. Suppose that $(\bar x, \bar y, \bar z, \bar u, \bar v, \bar w) \in \gph S^\nu \cap \ball(0, \rho)$, where the ball is specified by \eqref{eqn:normCompThm2}. Then,
\[
\bar u = F^\nu(\bar x) -\bar z, ~~~~ \bar v + \bar y \in \partial h^\nu(\bar z), ~~~~\bar w \in \nsum_{i=1}^m \bar y_i \con \partial f_i^\nu(\bar x) + N_{X}(\bar x).
\]
We construct a point $(\bar x, \bar y, z, u, v, w) \in \gph S$ that is near $(\bar x, \bar y, \bar z, \bar u, \bar v, \bar w)$ in the norm \eqref{eqn:normCompThm2}.

Since $(\bar x, \bar y, \bar z, \bar u, \bar v, \bar w) \in \ball(0, \rho)$, one has $\max\{\|\bar z\|_2, \|\bar v+\bar y\|_2\} \leq 2\rho$. Moreover, $\gph\partial h$ is nonempty, which follows largely from the definition of subgradients; cf. \cite[Corollary 8.10]{VaAn}. These facts ensure that there are $z,v\in\reals^m$ such that $(z,v+\bar y) \in \gph \partial h$ and
\[
\max\Big\{ \|z-\bar z\|_2, \big\|(v + \bar y) - (\bar v+\bar y)\big\|_2\Big\} \leq \exs_{2\rho}(\gph \partial h^\nu;~ \gph \partial h).
\]
There are $\bar a_i \in \con \partial f_i^\nu(\bar x)$ such that $\bar w - \nsum_{i=1}^m \bar y_i \bar a_i \in N_X(\bar x)$. By Caratheodory's theorem, there exist $\lambda_i^k\geq 0$ and $a_i^k \in \partial f_i^\nu(\bar x)$, $k=0, 1, \dots, n$, such that $\sum_{k=0}^n \lambda_i^k = 1$ and $\bar a_i = \sum_{k=0}^n \lambda_i^k a_i^k$. By assumption, there exists $b_i^k \in \con \partial f_i(\bar x)$ such that $\|a_i^k - b_i^k\|_2 \leq \eta^\nu$ and this holds for each $i$ and $k$. Let $\bar b_i = \sum_{k=0}^n \lambda_i^k b_i^k$, which then also is in $\con \partial f_i(\bar x)$. We construct
\[
w = \bar w + \nsum_{i=1}^m \bar y_i (\bar b_i - \bar a_i) ~~~\mbox{ and } ~~~ u = F(\bar x) - z.
\]
It now follows that $(\bar x, \bar y, z, u, v, w) \in \gph S$ and it remains to compute the distance between this point and $(\bar x, \bar y, \bar z, \bar u, \bar v, \bar w)$. We already have $\max\{\|z-\bar z\|_2, \|v - \bar v\|_2\}  \leq  \exs_{2\rho}(\gph \partial h^\nu;~ \gph \partial h)$. Moreover,
\[
\|u-\bar u\|_2  \leq \big\|F^\nu(\bar x) - F(\bar x)\big\|_2 + \|z - \bar z\|_2 \leq  \eta_0^\nu + \exs_{2\rho}(\gph \partial h^\nu;~ \gph \partial h).
\]
Finally, one has $\|w - \bar w\|_2 \leq \nsum_{i=1}^m |\bar y_i| \|\bar b_i - \bar a_i\|_2 \leq \sqrt{m} \rho \eta^\nu$. \eop

The theorem allows us to bound the solution error associated with a particular approximation in terms of the error in the individual components $h^\nu$ and $F^\nu$. (We refer to \cite{Royset.20b} for error estimates involving $X^\nu \neq X$ in the simplified setting with smooth $F$.) Thus, the focus turns to estimating the component errors. For that purpose, the following fact, which is a direct application of Theorem 5.2 in \cite{AttouchWets.91}, is useful. It leverages the  {\em truncated Hausdorff distance} between the sets $C$ and $D$:
\[
\hatsetd_\rho(C,D) = \max\big\{ \exs_\rho(C;D);\,\exs_\rho(D;C)\big\}.
\]

\begin{proposition}{\rm (approximation of subgradients).}\label{psubgradmappings}
Suppose that the norm on $\reals^{m+1}$ is $\max\{\|z\|_2,|\alpha|\}$, the norm on $\reals^m\times\reals^m$ is $\max\{\|z\|_2, \|v\|_2\}$, and  $\rho>\max\{ \dist(0,\epi h^\nu), \dist(0,\epi h)\}$. Then, there are $\alpha,\bar\rho\in [0,\infty)$, dependent on $\rho$, such that
\[
\exs_\rho\big(\gph \partial h^\nu; \,\gph \partial h\big) \leq \hatsetd_\rho\big(\gph \partial h^\nu, \,\gph \partial h\big) \leq \alpha \sqrt{\hatsetd_{\bar \rho}(\epi h^\nu,\epi h)} \leq \alpha \sqrt{\nsup_{z\in \ball(0,\bar\rho)}\big|h(z) - h^\nu(z)\big|}.
\]
\end{proposition}

The truncated Hausdorff distance between epigraphs can be computed using results from \cite{Royset.20b}, with the sup-expression given above being among the crudest possibilities.

\begin{example}{\rm (goal optimization; cont.).}\label{eGoal2}
In Example \ref{eGoal}, smoothing of $h^\nu$ causes a solution error
\[
\exs_\rho\Big(\,(S^\nu)^{-1}\big(\ball(0, \delta^\nu)\big); \,~S^{-1}\big(\ball(0, \epsilon)\big)    \,\Big) \leq \beta /\sqrt{\theta^\nu}
\]
for a constant $\beta$ when $0 \leq \delta^\nu \leq \rho<\infty$, $\epsilon \geq  \delta^\nu + \beta/\sqrt{\theta^\nu}$, and $\rho > \half \max\{\sqrt{\sum_{i=1}^m \min\{0, \tau_i\}^2}, \bar\alpha/\theta^\nu\}$, where $\bar\alpha = (\ln 2)(\sum_{i=1}^m \alpha_i)$.
\end{example}
\state Detail. Since $\nsup_{z\in \reals^m} |h^\nu(z) - h(z)| \leq \bar\alpha/ \theta^\nu$ by Example \ref{eGoal}, we obtain via Proposition \ref{psubgradmappings} that
\[
\exs_{2\rho}(\gph \partial h^\nu;\, \gph \partial h) \leq \alpha\sqrt{\bar\alpha/ \theta^\nu}
\]
for some $\alpha\in [0,\infty)$ when $\rho$ is sufficiently large as indicated. The claim then follows by Theorem \ref{tConstentApproxCompQuant} and Proposition \ref{thm:approxgeneralizedequations}.\eop

%
%

\begin{example}{\rm (distributionally robust optimization; cont.).}\label{eDistr2}
For Example \ref{eDistr}, there is a constant $\beta\in [0,\infty)$ such that
\[
\exs_\rho\Big(\,(S^\nu)^{-1}\big(\ball(0, \delta^\nu)\big); \,~S^{-1}\big(\ball(0, \epsilon)\big)    \,\Big) \leq \beta\sqrt{\alpha^\nu}
\]
as long as $0 \leq \delta^\nu \leq \rho<\infty$ and $\epsilon \geq \delta^\nu + \beta\sqrt{\alpha^\nu}$, where $\alpha^\nu = \max\{\exs(A; A^\nu), \exs(A^\nu; A)\}$.
\end{example}
\state Detail. Since $\dist(0, \epi h^\nu) = \dist(0, \epi h) = 0$, it follows by Proposition \ref{psubgradmappings} that
\begin{equation}\label{eqn:robustderiv1}
\hatsetd_{2\rho}(\gph \partial h^\nu, \gph \partial h) \leq \alpha \sqrt{\nsup_{z\in \ball(0,\bar\rho)}\big|h(z) - h^\nu(z)\big|}
\end{equation}
for some $\alpha,\bar \rho\in [0,\infty)$. Fix $z\in\ball(0,\bar\rho)$ and let $p^\nu \in \hat A^\nu(z)$; see Example \ref{eDistr} for notation. Then, there exists $\bar p\in A$ such that $\|p^\nu - \bar p \|_2 \leq \alpha^\nu$. Thus, $h^\nu(z) - h(z) \leq \langle p^\nu -\bar p, z\rangle \leq \|p^\nu -\bar p\|_2\|z\|_2 \leq \bar\rho \alpha^\nu$. We repeat this argument with the roles of $h$ and $h^\nu$ reversed and conclude that $|h(z) - h^\nu(z)| \leq \bar\rho \alpha^\nu$. In combination with \eqref{eqn:robustderiv1}, Theorem \ref{tConstentApproxCompQuant}, and Proposition \ref{thm:approxgeneralizedequations}, the claim follows.\eop

\begin{example}{\rm (augmented Lagrangian methods; cont.).}\label{eAug2}
For Example \ref{eAug}, we obtain that
\[
\exs_\rho\Big(\,(S^\nu)^{-1}\big(\ball(0, \delta^\nu)\big); \,~S^{-1}\big(\ball(0, \epsilon)\big)    \,\Big) \leq \beta^\nu/\theta^\nu, ~\mbox{ where } \beta^\nu = \big(2\rho + \|y^\nu\|_\infty\big)\sqrt{m-1},
\]
provided that $0 \leq \delta^\nu \leq \rho <\infty$ and $\epsilon \geq \delta^\nu + \beta^\nu/\theta^\nu$.
\end{example}
\state Detail. In this case, $\gph \partial h = ( \reals \times \{0\}^{m-1} ) \times (\{1\} \times \reals^{m-1} )$ and \[
\gph \partial h^\nu = \big\{(z,v) \in \reals^{2m}~\big|~v = (1, y_2^\nu + \theta^\nu z_2, \dots,  y_m^\nu + \theta^\nu z_m)\big\}.
\]
In view Proposition \ref{thm:approxgeneralizedequations} and Theorem \ref{tConstentApproxCompQuant}, it suffices to bound $\exs_{2\rho}(\gph \partial h^\nu; \gph \partial h)$. Let $(z^\nu, v^\nu)\in \gph \partial h^\nu \cap \ball(0,2\rho)$. We construct $(\bar z, \bar v) \in \gph \partial h$ by setting $\bar v_1 = 1$, $\bar z_1 = z_1^\nu$, $\bar z_i = 0$,  $\bar v_i = y_i^\nu + \theta^\nu z_i^\nu$ for $i=2, \dots, m$. Since $\|v^\nu\|_2\leq 2\rho$, one has $|v_i^\nu|\leq 2\rho$  and $ |y_i^\nu + \theta^\nu z_i^\nu| \leq 2\rho$ for $i=2, \dots, m$. Moreover, $\|v^\nu - \bar v\|_2^2 = 0$ and
\[
\|z^\nu - \bar z\|_2^2 = \nsum_{i=2}^m (z_i^\nu - \bar z_i)^2 = \nsum_{i=2}^m (z_i^\nu)^2 \leq \nsum_{i=2}^m\Big(\frac{2\rho + |y_i^\nu|}{\theta^\nu}\Big)^2 \leq (\beta^\nu/\theta^\nu)^2.
\]
Thus, the distance (in the appropriate norm) between  $(z^\nu,v^\nu)$ and $(\bar z, \bar v)$ is at most $\beta^\nu/\theta^\nu$.\eop

\begin{example}{\rm (exact penalty methods).}\label{eExact}
We next approach the problem in Example \ref{eAug} using an exact penalty method. Thus, $h(z) = z_1 + \sum_{i=2}^m \iota_{\{0\}}(z_i)$ as before, but for $\theta^\nu \in [0,\infty)$ we set
\[
h^\nu(z) = z_1 + \nsum_{i=2}^m \theta^\nu |z_i|.
\]
If $\theta^\nu \to \infty$, then the resulting approximations are consistent; the details are omitted. Interestingly,
\[
(S^\nu)^{-1}\big(\ball(0, \delta^\nu)\big) \cap\ball(0,\rho) \subset S^{-1}\big(\ball(0, \epsilon)\big)
\]
provided that $0 \leq \delta^\nu \leq \rho <\infty$, $\epsilon \geq \delta^\nu$, and $\theta^\nu \geq 2\rho$. This means that the approximating problem has a certain ``exactness'' property: for sufficiently large $\theta^\nu$, satisfaction of the optimality condition for the approximation problem locally in $\ball(0,\rho)$ implies satisfaction of the condition for the actual problem.
\end{example}
\state Detail. Example \ref{eAug2} furnishes $\gph \partial h$. Moreover, $\partial h^\nu(z) = \{1\} \times C_2 \times \cdots \times C_m$, where $C_i = \{\theta^\nu\}$ if $z_i > 0$, $C_i = [-\theta^\nu, \theta^\nu]$  if $z_i = 0$, and $C_i = \{-\theta^\nu\}$ otherwise.

Suppose that $\theta^\nu \geq 2\rho$. Let $(z^\nu,v^\nu)\in \gph\partial h^\nu \cap \ball(0,2\rho)$. Then, $v_1^\nu = 1$ and $v_i^\nu \in C_i$ for $i = 2, \dots, m$. Since $\|v^\nu\|_2 \leq 2\rho$, one has $|v_i^\nu| < 2\rho\leq \theta^\nu$ for $i=2, \dots, m$. This means that $z_i^\nu = 0$ for such $i$ in view of the definition of $C_i$. We conclude that $(z^\nu,v^\nu) \in \gph \partial h$ and $\exs_{2\rho}(\gph \partial h^\nu;\, \gph \partial h) = 0$. This in turn implies that $\exs_\rho(\gph S^\nu; \, \gph S) = 0$ by Theorem \ref{tConstentApproxCompQuant}. The claim follows from Proposition \ref{thm:approxgeneralizedequations}.\eop

\begin{example}{\rm (homotopy method; cont.).}\label{eHomotopy2}
In Example \ref{eHomotopy}, we obtain for $\lambda^\nu \in (0,1)$ that
\[
\exs_\rho\Big(\,(S^\nu)^{-1}\big(\ball(0, \delta^\nu)\big); \,~S^{-1}\big(\ball(0, \epsilon)\big)    \,\Big) \leq \beta^\nu \lambda^\nu, \mbox{ where } \beta^\nu = \sqrt{1+\frac{4\rho^2-(\lambda^\nu)^2}{(1-\lambda^\nu)^2}},
\]
provided that $0 \leq \delta^\nu \leq \rho <\infty$, $\epsilon \geq \delta^\nu + \beta^\nu\lambda^\nu$, and $\rho\geq \lambda^\nu/2$.
\end{example}
\state Detail. In this case, $\partial h(z) = \partial \hat h(z_1, \dots, z_{m-1}) \times \{0\}$ and $\partial h^\nu(z) = (1-\lambda^\nu)\partial \hat h(z_1, \dots, z_{m-1}) \times
\{\lambda^\nu\}$. To bound $\exs_{2\rho}(\gph \partial h^\nu;\, \gph \partial h)$, let $(z^\nu,v^\nu) \in \gph \partial h^\nu \cap \ball(0,2\rho)$. Then $v^\nu =
((1-\lambda^\nu)\hat y^\nu, \lambda^\nu)$ with $\hat y^\nu\in \partial \hat h(z_1^\nu, \dots, z_{m-1}^\nu)$. We construct $(\bar z,\bar v) \in \gph \partial h$ by setting $\bar v = (\hat y^\nu, 0)$ and $\bar z = z^\nu$. Then, one has $\|z^\nu - \bar z\|_2 = 0$ and $\|v^\nu - \bar v\|_2^2 = (1+\|\hat y^\nu\|_2^2)(\lambda^\nu)^2$. Since $\|v^\nu\|_2^2 = (1-\lambda^\nu)^2 \|\hat y^\nu\|_2^2 + (\lambda^\nu)^2 \leq 4\rho^2$, we obtain that $1+\|\hat y^\nu\|_2^2 \leq (\beta^\nu)^2$. From this we conclude that $\exs_{2\rho}(\gph \partial h^\nu;\, \gph \partial h) \leq \beta^\nu\lambda^\nu$. The result then follows by Theorem \ref{tConstentApproxCompQuant} and Proposition \ref{thm:approxgeneralizedequations}.\eop

\state Acknowledgement. The author is thankful to K. Balasubramanian (UC Davis) for introducing him to the inverse problem in Section 5. The research is supported in part by ONR Science of Autonomy (N0001421WX00142) and AFOSR (18RT0599, 21RT0484).


\bibliographystyle{plain}
\bibliography{refs}

\end{document}